\documentclass[11pt]{article}
\usepackage{epsfig}
\usepackage{amsmath,amssymb, graphicx, enumitem, mathrsfs, mathtools, breqn, booktabs, relsize,hyperref,float,amsthm}
\usepackage[version=4]{mhchem}
\usepackage[authoryear]{natbib}
\usepackage{listings,wrapfig,tabu}
\usepackage[caption=false]{subfig}
\usepackage{color}
\usepackage{xcolor}
\usepackage{fullpage}
\title{Analysis and Simulation of a Novel Run-and-Tumble Model with Autochemotaxis}
\author{Nicholas J. Russell and Louis F. Rossi}
\date{February 3, 2021}

\newtheorem{proposition}{Proposition}

\DeclareMathOperator*{\argmax}{arg\,max}
\begin{document}

\maketitle
\begin{abstract}
We model, analyze, and simulate a novel run-and-tumble model with autochemotaxis, biologically inspired by the phytoplankton \textit{Heterosigma akashiwo}. Developing a fundamental understanding of planktonic movements and interactions through phototaxis and chemotaxis is vital to comprehending why harmful algal blooms (HABs) start to form and how they can be prevented. We develop a one- and two-dimensional mathematical and computational model reflecting the movement of an ecology of plankton, incorporating both run-and-tumble motion and autochemotaxis. We present a succession of complex and biologically meaningful models combined with a sequence of laboratory and computational experiments that inform the ideas underlying the model. By analyzing the dynamics and pattern formation which are similar to experimental observations, we identify parameters that are significant in plankton's pattern formation in the absence of bulk fluid flow. We find that the precise form of chemical deposition and plankton sensitivity to small chemical gradients are crucial parameters that drive nonlinear pattern formation in the plankton density.  
\end{abstract}
\section{Introduction}

\indent Plankton is a general term for a wide range of floating or drifting organisms in fresh and salt water. They are not strong enough to swim against currents. Plankton are of considerable scientific interest for a 
number of reasons. For one, they are an important food source for larger 
organisms. Also, they process about 60\% of the marine carbon, and they consume a significant fraction of man-made \ce{CO2} while producing about half of the atmospheric oxygen (\citealt{Berdalet2016,Riebesell2008}).
Nonetheless, plankton exhibit a rich array of behaviors including self-propulsion, aggregation, photosynthesis, feeding, and predation (\citealt{Sheng2007,Sheng2010}). 
Sometimes, coordinated microscopic behaviors lead to significant mesoscale processes. Most notably, several types of plankton are responsible for harmful algal blooms (HABs) wherein the population density of plankton skyrockets in a concentrated area (\citealt{Hallegraeff2004}). However, some types of plankton, and particularly those involved in the formation of HABs, aggregate even when there is no external current. While it is reasonable to model plankton as independent passive scalar quantities drifting in a flow field, a more detailed understanding of their behavior may help scientists understand and anticipate when plankton will aggregate and what combinations of features will drive this process.

In this paper, we seek to construct, analyze, and simulate a run-and-tumble model with autochemotaxis. This work is biologically inspired by the dynamics of the phytoplankton \textit{Heterosigma akashiwo}. Over the past several years, this raphidophyte has caused several HABs throughout the world, including Japan, Scandanavia, South Carolina, and the Delaware Bay (\citealt{Hennige2013}). These HABs are responsible for large fish kills (e.g. the Puget Sound in 2006) and public health hazards for the communities that surround the coasts (\citealt{Lewitus2012}). This is of particular interest to the Delaware community, as these HABs are becoming more frequent over the past several years in the Delaware Bay. \textit{Heterosigma akashiwo} are predatory swimmers, capable of consuming bacteria and performing an inefficient photosynthesis. Even in the absence of prey or bulk fluid motion, plankton exhibit complex pattern formations in two- and three-dimensions. Their behavior is a direct response of local chemical, turbulent, and optical cues, so a fundamental understanding of how they swim and hunt is essential to anticipating changes in marine ecology in response to climate change and ocean acidification (\citealt{Fu2012,Sengupta2017}).

We discuss a research program at the University of Delaware to mathematically model, analyze, and understand complex plankton behavior. These activities include careful experiments, isolating specific variables and behaviors to connect them to a unique coupling of well-established processes in predation, chemotaxis, and phototaxis. We develop and analyze a model for plankton motion and signaling that is built upon the works of several authors (\citealt{Hillen2002,Hillen2009a,Keller1971,Keller1971a,Lushi2012a,Bearon2000a}) in order to gain insight and understanding into the interacting mechanisms of organism propulsion and chemical deposition, diffusion, and decay.

The motion of plankton is driven by a number of processes, only some of which have been fully explored experimentally for our model organism, \textit{Heterosigma akashiwo}. Therefore, it is necessary to limit the study of their motion to a subset of these which have the greatest impact on their bulk coordinated behavior. In this paper, we will limit our discussion to run-and-tumble propulsive behavior which has been studied extensively (\citealt{Keller1971,Keller1971a}).
When combined with autochemotaxis, in which organisms release the chemical they are attracted to, run-and-tumble motion admits quite rich dynamics which mirror experimental findings. We develop and analyze a one- and two-dimensional model of plankton motion incorporating run-and-tumble motion along with autochemotaxis while using a biologically relevant tumbling probability. We acknowledge that past studies have suggested that gravity and rotational diffusion can describe the movement of \textit{Heterosigma akashiwo} in three-dimensions (e.g. \citealt{Bearon2008}). However, gravity and rotational diffusion are negligible in two dimensions, and movement can be described well by run-and-tumble motion (\citealt{Visser2008,Shaw2017}).

In Section \ref{sec:Experiments}, we summarize the experiments conducted and the resultant observations which inform our mathematical model. In Section \ref{sec:model}, we develop analytical models of run-and-tumble dynamics coupled with autochemotaxis in one and two dimensions. In Section \ref{sec:Analysis1D}, an analysis of the one-dimensional model is conducted both analytically and numerically. In Section \ref{sec:Analysis2D}, we analyze our two-dimensional model analytically and numerically. In Section \ref{sec:Discussion}, we discuss the analysis and simulations by connecting our results to experimental observations. Finally, in Section \ref{sec:Conclusions}, we share our conclusions of this work and future research opportunities.

\section{Some experimental insights into plankton aggregation}\label{sec:Experiments}
In this paper, we use \textit{Heterosigma akashiwo}, a constituent plankton in an algal bloom colloquially referred to as the ``red tide", as our model organism. To gain insight on the evolution and dynamics of the plankton distributions, we performed a series of laboratory experiments using a camera to gather raw data. See Appendix A for the full details on experimentation. The data should be considered qualitative in the sense that we could not rigorously map the image density to a calibrated population density.



We anticipate that two main mechanisms drive the plankton density distribution in water without external fluid flow: chemotaxis and phototaxis. For chemotaxis, plankton swim in response to a chemical while also emitting the chemical, also called autochemotaxis (\citealt{Seymour2009}). Phototaxis is dependent on the plankton's photosystem behavior and the depth of the water, and light is a source of energy for the organism. Its photosystem is degraded by light and is constantly repaired by internal biochemical reactions (\citealt{Fu2008}). As well, plankton are observed to vertically and horizontally migrate towards and away from light, and it is reasonable to hypothesize that the state of their photosystem drives this behavior (\citealt{Hara1987}). Lastly, recent studies have shown that plankton are able to change their migration strategies due to microscopic turbulent cues and thus, we should not model them as passive scalars in a fluid flow (\citealt{Sengupta2017}). It is our aim to investigate autochemotactic behavior in this model organism as a foundation for exploring its richer behavior in response to phototaxis and other processes impacting migration and aggregation.
\begin{figure}[t!]
    \centering
    \subfloat[Deep plate experiment]{
        \includegraphics[width=0.33\linewidth, angle=180,origin=c]{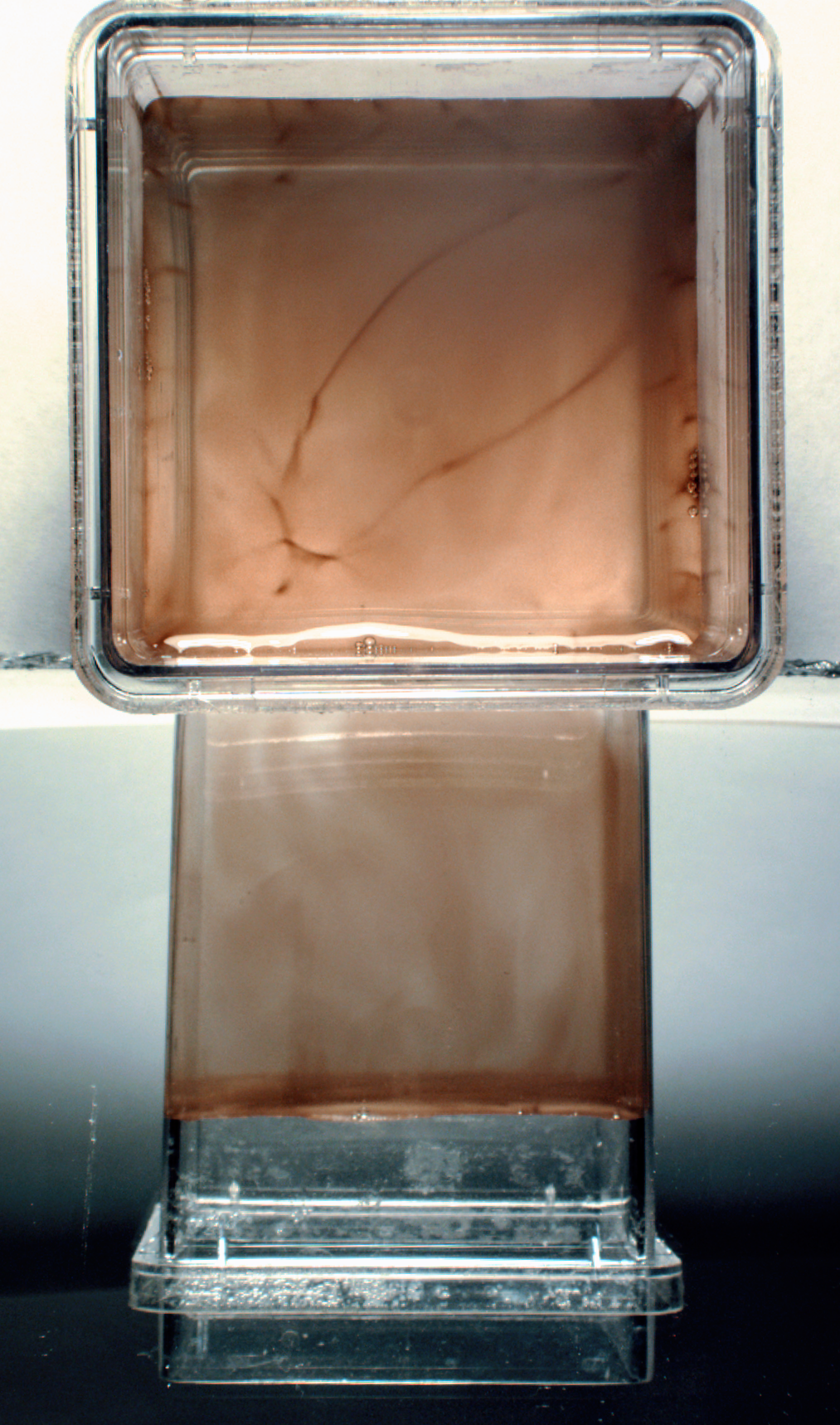}\label{fig:DeepPlate}}
        \hspace{0.1\linewidth}
    \subfloat[Shallow plate experiment]{
        \includegraphics[width=0.4\linewidth]{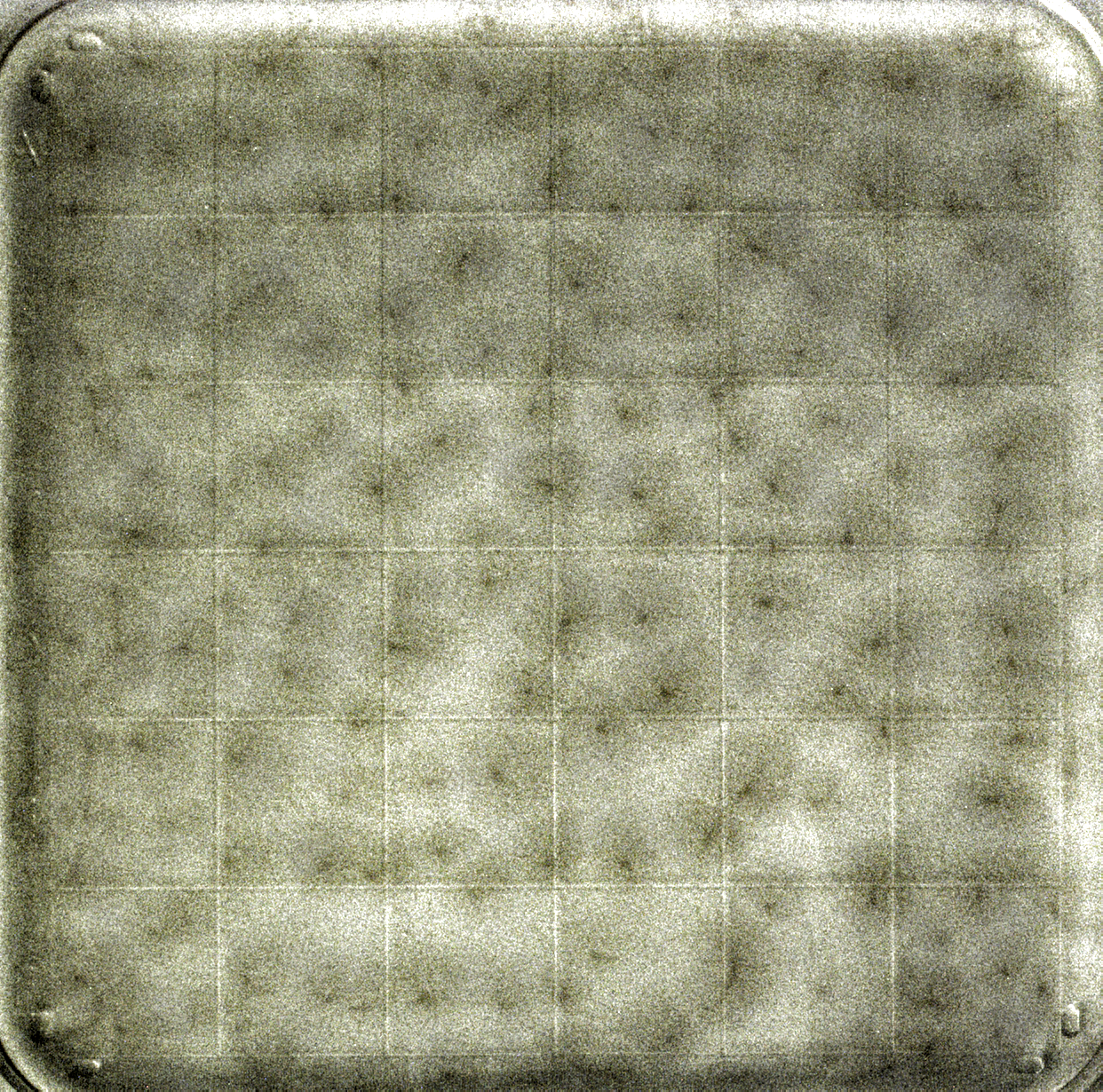}\label{fig:DarkExp}}
    \caption{Pictures from various experiments conducted. The aggregations in \textbf{b} formed after 5 minutes in complete darkness (an infrared camera was used).}
    \label{fig:DeepDark}
\end{figure}

To demonstrate the richness of the full dynamics of \textit{Heterosigma akashiwo}, even in the absence of an external driving fluid flow, we conducted experiments in a deep tank. In this deep tank, all effects from chemotaxis and phototaxis are observable.  We observe plankton aggregating into complex patches and filaments, while migrating towards and away from the light (see Fig. \ref{fig:DeepPlate}). Water attenuates light, captured by the Lambert-Beer Law, and aggregations can shade other plankton.

\begin{figure}
\centering
\subfloat[Time Elapsed: 2 Minutes]{\includegraphics[width=0.4\linewidth]{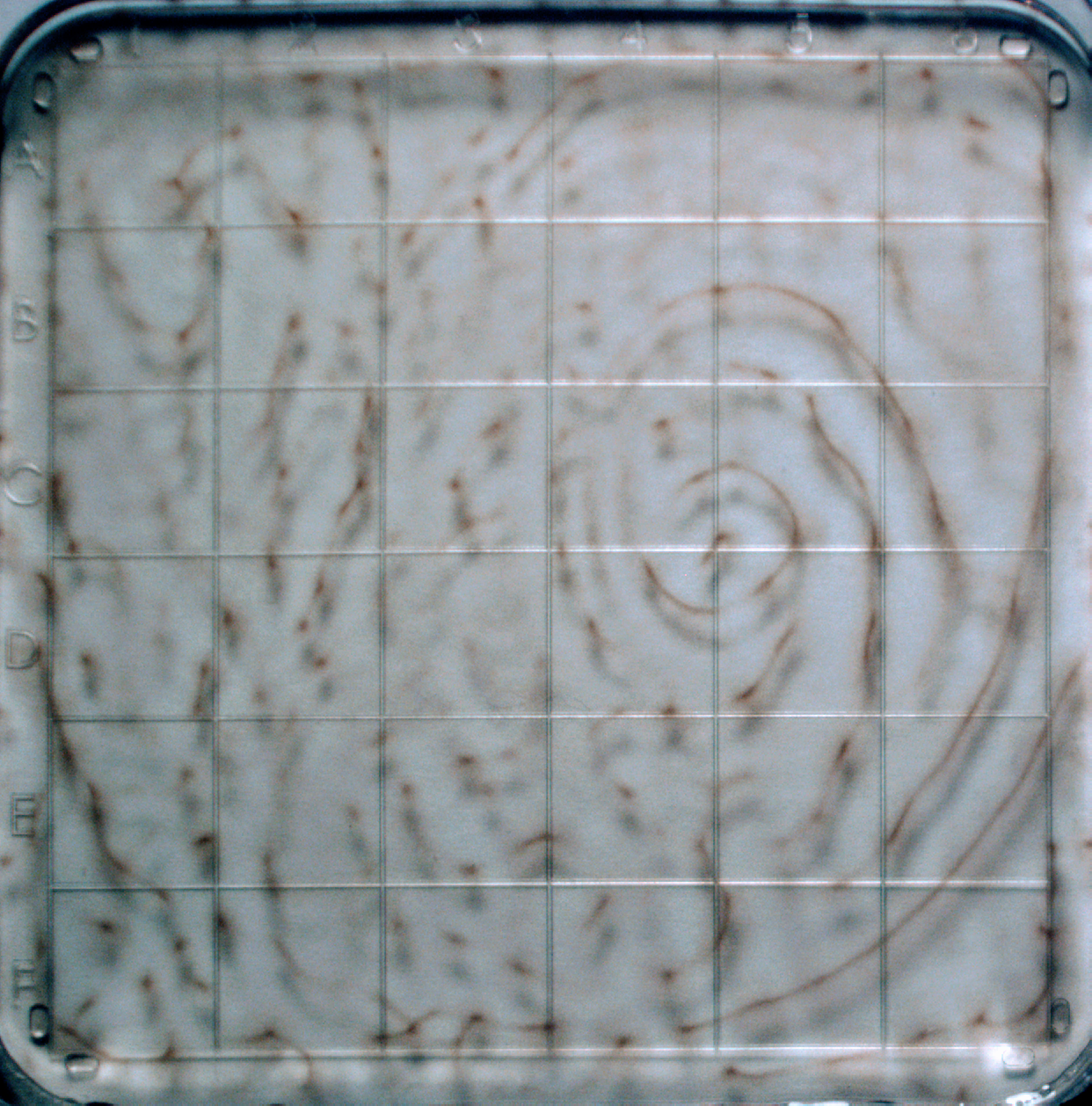}}\hspace{0.1\linewidth}
\subfloat[Time Elapsed: 6 Minutes]{\includegraphics[width=0.4\linewidth]{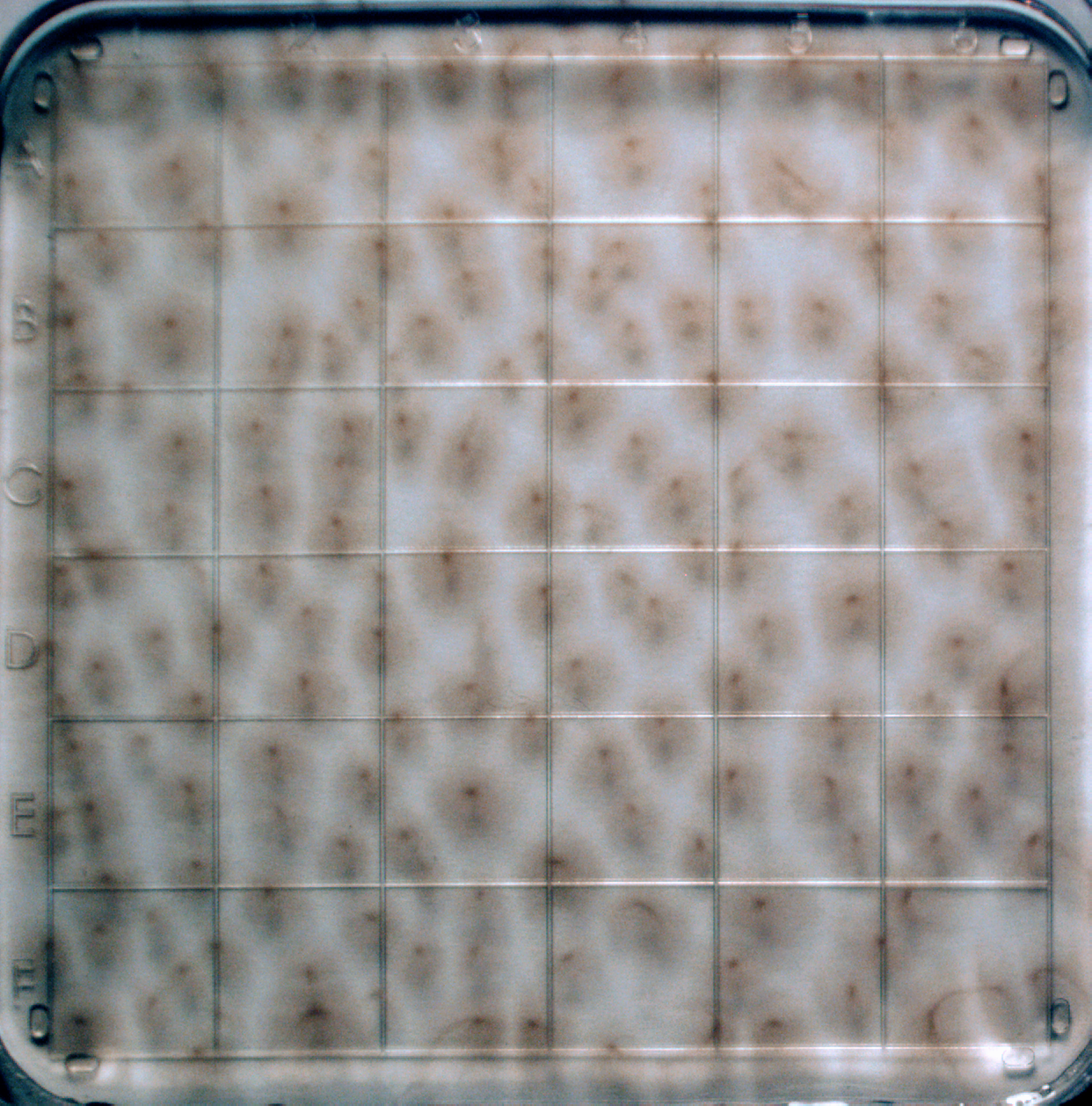}}\\
\subfloat[Time Elapsed: 15 Minutes]{\includegraphics[width=0.4\linewidth]{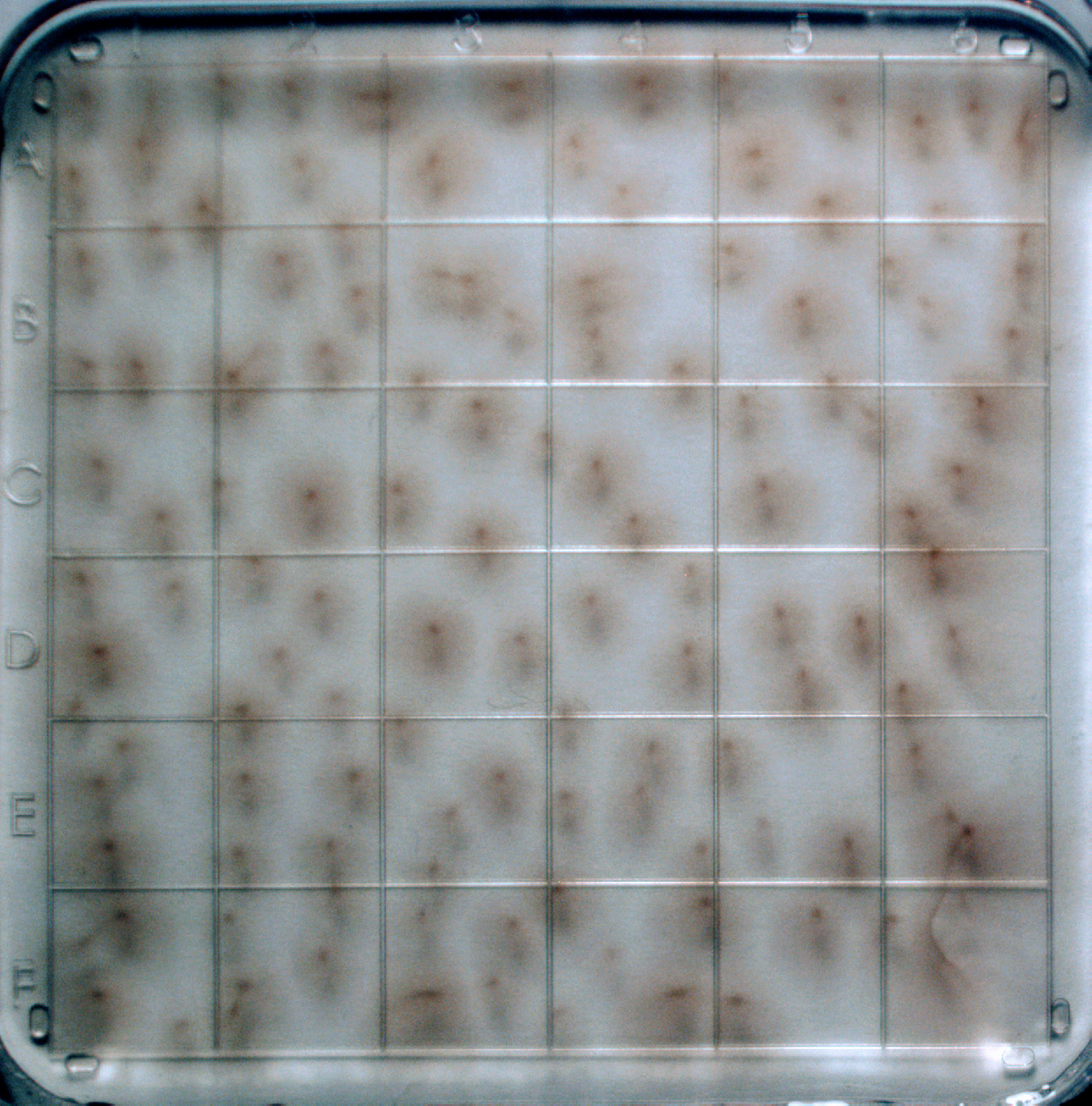}}\hspace{0.1\linewidth}
\subfloat[Time Elapsed: 40 Minutes]{\includegraphics[width=0.4\linewidth]{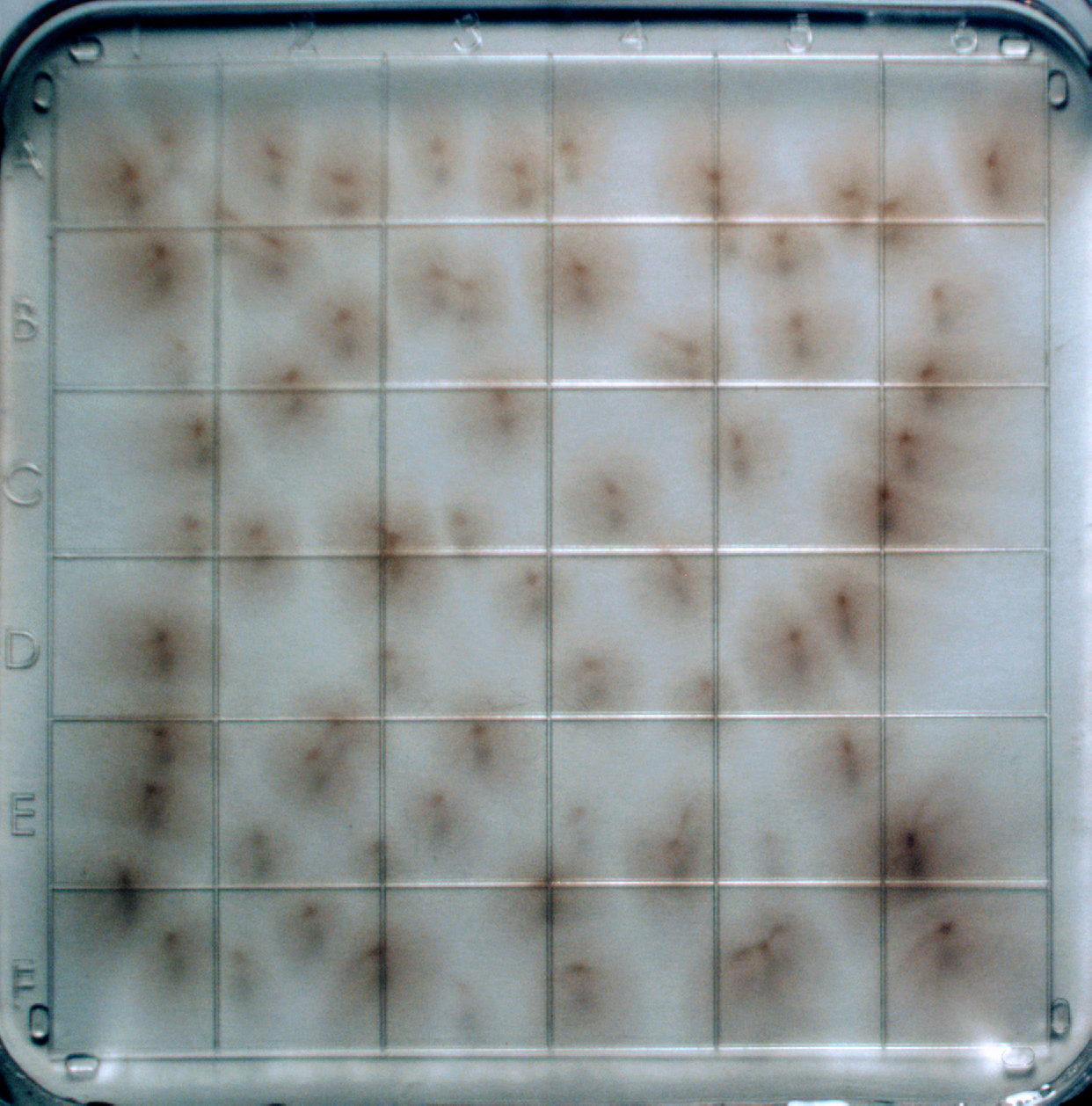}}\\
\subfloat[Time Elapsed: 60 Minutes]{\includegraphics[width=0.4\linewidth]{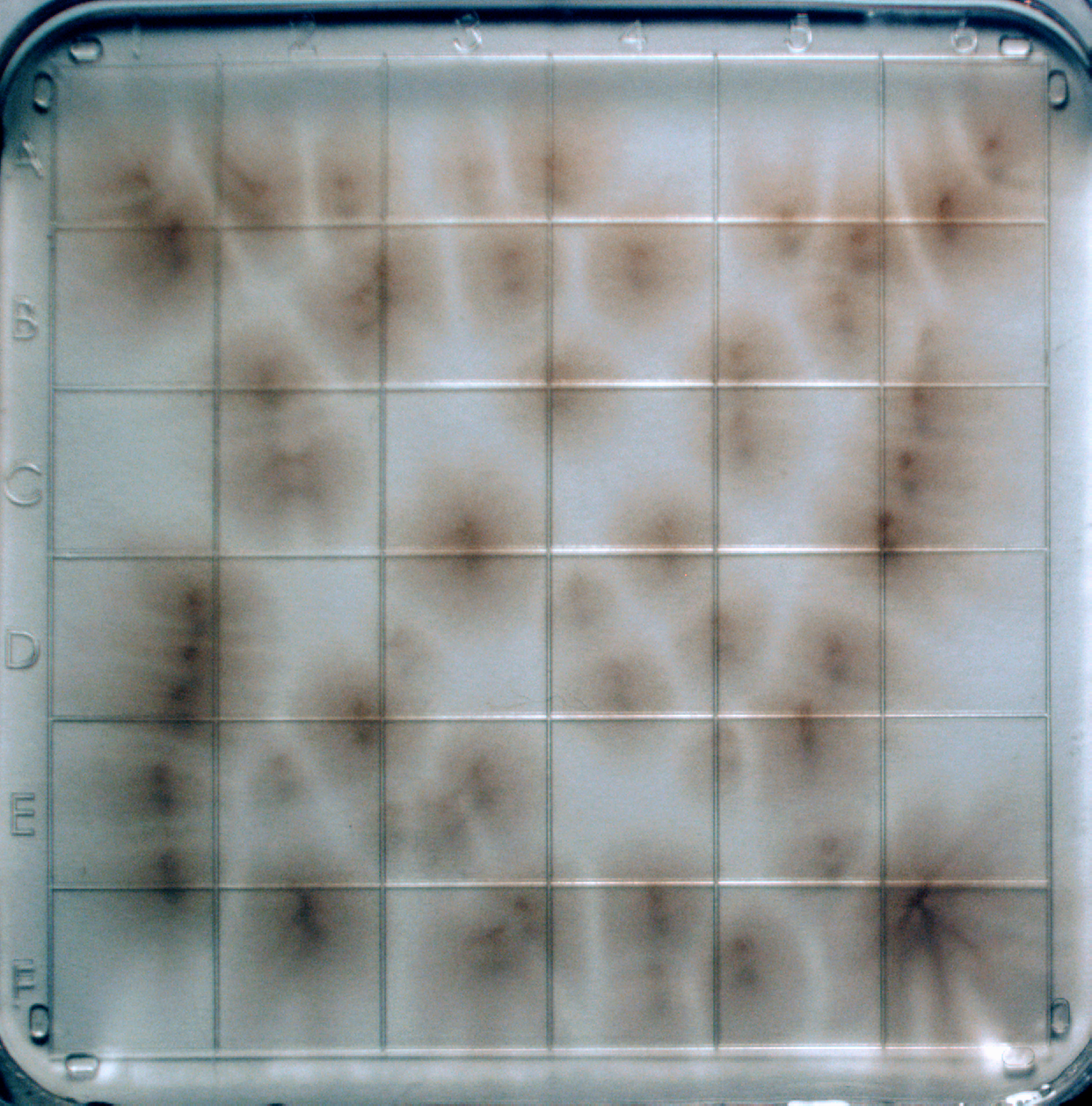}\label{fig:T60}}\hspace{0.1\linewidth}
\begin{minipage}[b]{0.4\linewidth}
    \caption{The evolution of aggregation patterns from \textit{Heterosigma akashiwo} during an hour long experiment. We used a square, shallow dish with a length of 6.3 cm along with a lamp placed at a $60^{\circ}$ angle of incidence. We used a filter on each picture to make the aggregations clearer. We note that the spiral pattern in \textbf{a} is caused by initial fluid motion when the plankton solution is poured into the dish initially.}
    \label{fig:Progression}
\end{minipage}
\end{figure}
To remove optical attenuation and shading from consideration and to minimize the impact of vertical migration, the depth of the water column is reduced until the vertical effects were no longer observable, which occurs when the depth is less than 1 cm. All experiments we discuss from now on only consider a shallow dish to remove vertical migration from consideration as a dominant affect. The system is flooded with light, and thus there is no photo-intensity gradient in the field.  Therefore, the shallow environment allows us to eliminate certain affects attributable to the organisms sensitivity to light. We still must account for the fact that \textit{Heterosigma akashiwo} responds to light, even if light is uniformly distributed in the shallow environment.

To further isolate the behavioral impact of photosynthesis, we conduct an experiment where we film \textit{Heterosigma akashiwo} in complete darkness, allowing the photosystems in all organisms to reset, and then observe what occurs once light returns. Within the first 10 minutes, aggregations form in complete darkness (see Fig. \ref{fig:DarkExp}). Thus, we can establish that autochemotaxis drives the formation of aggregations. Conducting experiments with lights on, we observe that phototaxis drives the coarsening and movement of these aggregations. From a simple experiment shown in Fig. \ref{fig:Progression}, aggregations started forming within 2 minutes, more rapidly than the time it takes the photosystems to degrade and break. The structure of the aggregations is also different with and without light, so phototaxis changes the dynamics of aggregation shapes.

We conducted experiments with a lid flush with the upper liquid surface to see if Marangoni effects were playing a role in the formation of aggregations and found no observed difference with all other free surface experiments. From these experiments and observations, we infer that although phototaxis stimulates motion and perhaps aggregation, chemotaxis combined with run-and-tumble is the initial driver of aggregations of \textit{Heterosigma akashiwo}.

\section{A 1D and 2D model of plankton motion}\label{sec:model}

\subsection{The two-dimensional model}

Consider a system of $N_p$ plankton as point sources, 
\begin{equation}\label{eq:rhoinit}
    \rho(\vec{x}, t) = \sum\limits_{i= 1}^{N_p} \delta(\vec{x} - \vec{x}_i(t)), 
\end{equation}
where $\vec{x}_i(t) = (x_1(t), x_2(t))$ is the position of the $i$th plankton and $\delta(\vec{x})$ is the Dirac delta function. We denote $c(\vec{x}, t)$ as the chemical concentration at a point $\vec{x}$. When considering a classic run-and-tumble model for plankton motion, literature recommends the flipping probability of
\begin{equation}\label{eq:probflip}
    \mathbb{P}\left[ \text{Tumble in } (t, t + \Delta t) \right] = \frac{\lambda_0}{2}\left(1 - \frac{\vec{v} \cdot \nabla c}{\| \vec{v} \cdot \nabla c\| }\right) \Delta t
\end{equation}
where $\vec{v}$ is the velocity of a plankton, $\lambda_0 \geq 0$, and $\nabla c$ is the chemical gradient (\citealt{Saintillan2018}). This model is problematic from a biological perspective. We see that $\frac{\vec{v} \cdot \nabla c}{\| \vec{v} \cdot \nabla c\|} = \text{sgn}\left( \| \vec{v} \cdot \nabla c\|\right)$ takes three distinct values: $-1, 0,$ and 1. This gives a drastic shift in behavior when plankton flips to being orthogonal to the gradient, implying an infinite sensitivity to both the environment and the organism's response to the chemical. Thus, we propose to create a transition to smooth out this hard shift by using the parameter $\delta > 0$ and converting \eqref{eq:probflip} into 
\begin{equation}\label{eq:NewProbflip}
       \mathbb{P}\left[ \text{Tumble in } (t, t + \Delta t) \right] = \frac{\lambda_0}{2}\left(1 - \frac{\vec{v} \cdot \nabla c}{\sqrt{\left(\vec{v} \cdot \nabla c\right)^2 + \delta^2}} \right) \Delta t.
\end{equation}
Biologically, $\delta$ captures the ability of an organism to sense weak gradients in its environment. As $\delta \to 0$, the organism can perfectly sense an infinitesimally small gradient, and we recover our original model \eqref{eq:probflip}. Larger values of $\delta$ describe a smooth transition in the organism's behavior which could capture a noisier response to weak signals when averaged over multiple instances of the same organism or multiple organisms in close proximity. We shall see later that $\delta$ has a strong impact on macroscopic pattern formation.

We assume the plankton move in a run-and-tumble motion with chemotaxis, i.e. the probability of a plankton tumbling and changing directions is given by \eqref{eq:NewProbflip}. If a plankton does tumble, the probability that a plankton moves from a direction $\vec{v}$ to $\vec{v}^{\prime}$ is defined by the transition probability $T(\vec{v}^{\prime}, \vec{v})$. We assume that plankton move at a constant speed $v$ and thus, $\vec{v} = v \vec{e}_{\theta}$, where $\vec{e}_{\theta} = \left\langle \cos(\theta), \sin(\theta) \right\rangle$. Because of this assumption, the transition probability can be rewritten as $T(\vec{v}^{\prime}, \vec{v}) = T(\theta^{\prime}, \theta)$. 

\indent The chemical attractant, which is produced by each plankton, diffuses and decays throughout the system. Since there is no consensus on when the plankton emit this chemical, we will assume that plankton emit this chemical based on how much chemical they can sense at their current position. Let $f(c)$ be a general deposition function that is only dependent on the chemical, $c$. The diffusion and decay of $c$ can be captured in an evolution equation for the chemical attractant: 
\begin{dmath}\label{eq:Chemical}
    c_t = \kappa \Delta c - \beta c + f(c) \rho,
\end{dmath}
where $\kappa \geq 0$ is the diffusion coefficient and $\beta \geq 0$ is the decay coefficient. 

In order develop a model for $\rho$ from an Eulerian perspective, we now define $\rho(\vec{x}, t)$ to be the density of the plankton at a given time and location in $\mathbb{R}^2$. We introduce $\psi(\vec{x}, \theta, t)$ as the density of plankton at time $t$ and position $\vec{x}$, moving in the $\theta$ direction. Thus, the evolution equation for $\psi$ can be written as a Fokker-Plank equation:
\begin{dgroup*}
\begin{dmath}\label{eq:RhoPsi}
\rho(\vec{x}, t) = \int_0^{2\pi} \psi(\vec{x}, \theta, t) \, d\theta
\end{dmath}
\begin{dmath}\label{eq:2ddensity}
\psi_t = - v \vec{e}_{\theta} \cdot \nabla \psi - \frac{\lambda_0}{2} \left[ 1 - \frac{v \left(\nabla c \cdot \vec{e}_{\theta}\right)}{\sqrt{ v^2 \left(\nabla c \cdot \vec{e}_{\theta}\right)^2 + \delta^2 }} \right] \psi + \frac{\lambda_0}{2} \int_0^{2\pi} \left[ 1 - \frac{v \left(\nabla c \cdot \vec{e}_{\theta^{\prime}}\right)}{\sqrt{ v^2 \left(\nabla c \cdot \vec{e}_{\theta^{\prime}}\right)^2 + \delta^2 }} \right] T(\theta, \theta^{\prime}) \psi(\vec{x}, \theta^{\prime}, t) \, d \theta^{\prime}.
\end{dmath}
\end{dgroup*}
In $\eqref{eq:2ddensity}$, the first term represents plankton that have moved away from position $\vec{x}$; the second term represents the plankton at $\vec{x}$ who tumbled to a different direction than $\theta$; and the third term represents the plankton that have tumbled at $\vec{x}$ to the direction of $\theta$ from $\theta^{\prime}$. Therefore, \eqref{eq:Chemical}-\eqref{eq:2ddensity} is the continuum two-dimensional run-and-tumble with autochemotaxis system. To non-dimensionalize this model, we denote the natural time-scale of $\tau = \frac{1}{\lambda_0}$ and the natural length scale of $L = \frac{v}{\lambda_0}$. Define the non-dimensional time and length as $t = \tau t^{\star}$, $\vec{x} = L \vec{x}^{\star}$, where the $\star$-notation denotes a non-dimensional parameter. To simplify further, we assume that the turning kernel is a uniform distribution, i.e. $T(\theta^{\prime}, \theta) = \frac{1}{2 \pi}$. Removing the star notation, we obtain the following non-dimensional model of autochemotactic-driven motion in 2D: 
\begin{dgroup}
\begin{dmath}\label{eq:2nond1}
    c_t = d_1 \Delta c - d_2 c + f(c) \rho
\end{dmath}
\begin{dmath}\label{eq:2nond2}
\rho(\vec{x}, t) = \int_0^{2\pi} \psi(\vec{x}, \theta, t) \, d\theta
\end{dmath}
\begin{dmath}\label{eq:2nond3}
\psi_t = - \vec{e}_{\theta} \cdot \nabla \psi - \frac{1}{2} \left[ 1 - \frac{ \nabla c \cdot \vec{e}_{\theta}}{\sqrt{\left(\nabla c \cdot \vec{e}_{\theta}\right)^2 + \delta^2 }}\right]  \psi + \frac{1}{4\pi} \int_0^{2\pi} \left[ 1 - \frac{ \nabla c \cdot \vec{e}_{\theta^{\prime}}}{\sqrt{\left(\nabla c \cdot \vec{e}_{\theta^{\prime}}\right)^2 + \delta^2 }}\right] \psi(\vec{x}, \theta^{\prime}, t) \, d \theta^{\prime},
\end{dmath}
\end{dgroup}
where $d_1 = \frac{\kappa \lambda_0}{v^2}$, $d_2 = \frac{\beta}{\lambda_0}$, and $f(c)$ has been scaled by $\frac{1}{\lambda_0}$. Throughout this paper, we will analyze and simulate the system in a large, square domain $[0, \ell] \times [0, \ell]$ with periodic boundary conditions to express the near limitless domain in which plankton truly live. Computational limitations necessarily require that we make $\ell$ finite. 

\subsection{The 1D model}
For a one-dimensional model, $\theta$ is limited to values of 0 (right-moving) and $\pi$ (left-moving), and $\psi$, $\rho$, and $c$ vary with $x \in \Omega \subset \mathbb{R}$ and $t \geq 0$ only. We can remove the $\theta$ dependence entirely by replacing $\psi(\theta, x, t)$ with the distinct functions $\psi^+(x,t)$ and $\psi^-(x,t)$, representing right- and left-moving plankton respectively. Equation \eqref{eq:NewProbflip} can be rewritten as 
\begin{equation}
    \mathbb{P}\left[ \text{Tumble in } (t, t + \Delta t) \right] = \frac{\lambda_0}{2} \left( 1 - \frac{c_x}{\sqrt{c_x^2 + \delta^2}} \right),
\end{equation}  
and we assume that the plankton move at a constant speed $v$. For the autochemotaxis part of the model, let $f(c)$ denote an arbitrary deposition function. We can then write the evolution equations as follows: 
\begin{dgroup}
\begin{dmath}[label={eq:chemical}]
c_t = \kappa c_{xx} - \beta c + f(c)(\psi^+ + \psi^-)
\end{dmath}
\begin{dmath}[label={eq:plus}]
\psi_t^+ = - v \psi_x^+ - \frac{\lambda_0}{2} \left( 1- \frac{c_x}{\sqrt{c_x^2 + \delta^2}} \right) \psi^+ + \frac{\lambda_0}{2} \left( 1 + \frac{c_x}{\sqrt{c_x^2 + \delta^2}} \right) \psi^-
\end{dmath}
\begin{dmath}[label={eq:minus}]
\psi_t^- = v \psi_x^- + \frac{\lambda_0}{2} \left( 1- \frac{c_x}{\sqrt{c_x^2 + \delta^2}} \right) \psi^+ - \frac{\lambda_0}{2} \left( 1 + \frac{c_x}{\sqrt{c_x^2 + \delta^2}} \right) \psi^-. 
\end{dmath}
\end{dgroup} 
Non-dimensionalizing, we denote the natural time-scale of $\tau = \frac{1}{\lambda_0}$ and the natural length scale of $L = \frac{v}{\lambda_0}$. Define the non-dimensional time and length as $t = \tau t^{\star}$ and $x = L x^{\star}$, where the $\star$-notation denotes a non-dimensional parameter. Removing the star notation, we obtain the following non-dimensional model of $\eqref{eq:chemical}-\eqref{eq:minus}$:
\begin{dgroup}
\begin{dmath}\label{eq:nd1}
c_t = d_1 c_{xx} - d_2 c + f(c) (\psi^+ + \psi^-)
\end{dmath}
\begin{dmath}\label{eq:nd2}
\psi_t^+ = - \psi_x^+ - \frac{1}{2} \left( 1- \frac{c_x}{\sqrt{c_x^2 + \delta^2}} \right) \psi^+ + \frac{1}{2} \left( 1 +\frac{c_x}{\sqrt{c_x^2 + \delta^2}} \right) \psi^- 
\end{dmath}
\begin{dmath}\label{eq:nd3}
\psi_t^- = \psi_x^- + \frac{1}{2} \left( 1- \frac{c_x}{\sqrt{c_x^2 + \delta^2}} \right) \psi^+ - \frac{1}{2} \left( 1 + \frac{c_x}{\sqrt{c_x^2 + \delta^2}} \right) \psi^-
\end{dmath}
\end{dgroup}
where  $d_1 = \frac{\kappa \lambda_0}{v^2}$ and $d_2 = \frac{\beta}{\lambda_0}$, and $f(c)$, which is scaled by $\frac{1}{\lambda_0}$, are all non-dimensional. If we define $\rho(x, t) = \psi^+ + \psi^-$ to be the total density, then manipulating $\eqref{eq:nd2}$ and $\eqref{eq:nd3}$ allows us to obtain a coupled PDE system for $\rho$ and $c$:
\begin{dgroup}
\begin{dmath}\label{eq:nd11}
c_t = d_1 c_{xx} - d_2 c + f(c) \rho 
\end{dmath}
\begin{dmath}\label{eq:nd12}
\rho_{tt} + \rho_t = \rho_{xx} - \frac{\partial}{\partial x} \left[ \frac{c_x}{\sqrt{c_x^2 + \delta^2}} \rho \right] 
\end{dmath}
\end{dgroup}
For both \eqref{eq:nd1}-\eqref{eq:nd3} and \eqref{eq:nd11}-\eqref{eq:nd12}, we will assume periodic boundary conditions and consider the domain $\Omega = [0,\ell]$, where $\ell > 0$. 
\begin{figure}[t!]
    \centering
    \includegraphics[height=0.4\linewidth]{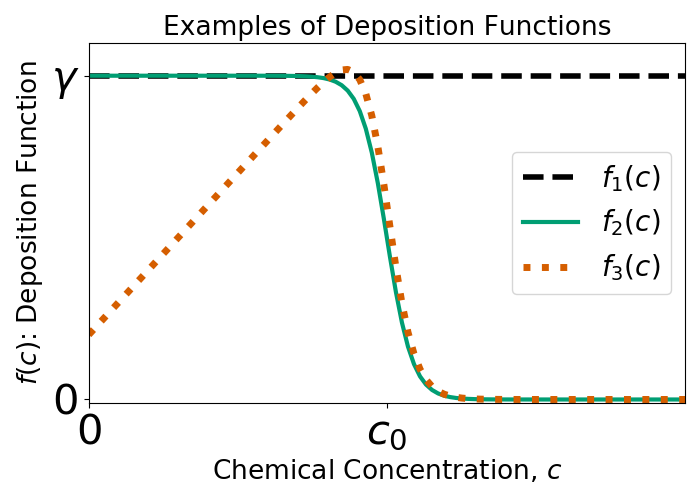}
    \caption{Examples of the three types of deposition functions considered in the simulations. The black dashed line is $f_1(c) = \gamma$, the constant deposition function; the green, solid line is $f_2(c)$, the switch deposition function; and the orange, dotted line is $f_3(c)$, the linear switch deposition functions. The analytic descriptions of $f_2(c)$ and $f_3(c)$ can be found in equations \eqref{eq:Switch} and \eqref{eq:LinearSwitch}, respectively.}
    \label{fig:Deposition}
\end{figure}
Since the deposition function for \textit{Heterosigma akashiwo} is unknown, we will utilize basic deposition functions related to other organisms such as insects and bacteria. The three deposition functions we analyze are $f_1(c) = \gamma$, $f_2(c) = \gamma \Theta(c_0 - c)$, and $f_3(c) =\gamma(c + \eta)\Theta(c_0 - c)$, where $\gamma, c_0, \eta > 0$ and $\Theta(x)$ is the Heaviside function. The constant function, $f_1$, is where an organism always emits the chemical; the switch, $f_2$, is where an organism releases the chemical only if the chemical concentration at $\vec{x}$ is below the threshold $c_0$; and the linear switch, $f_3$, is where an organism releases a specific amount of chemical dependent on the amount of chemical at $\vec{x}$, but does not emit if the chemical at $\vec{x}$ is above $c_0$. In our numerical simulations, we smooth out the Heaviside functions to avoid discontinuities as 
\begin{dgroup}
\begin{dmath}\label{eq:Switch}
f_2(c) = \frac{\gamma}{2} \left[ \tanh\left(\frac{c_0-c}{T_0} \right) + 1 \right]
\end{dmath}
\begin{dmath}\label{eq:LinearSwitch}
    \hspace{2mm} f_3(c) = \frac{\gamma(c + 0.2c_0)}{2 c_0} \left[ \tanh\left(\frac{c_0-c}{T_0} \right) + 1 \right],
\end{dmath}
\end{dgroup}
where $T_0 > 0$. Examples of these deposition function are displayed in Fig. \ref{fig:Deposition}. For a summary of all dimensional and non-dimensional parameters, see Tables \ref{tab:ParametersD} and \ref{tab:ParametersND}, respectively.
\begin{table}[t!]
\centering
\renewcommand\arraystretch{1.3}
\begin{tabular}{|c|c|c|c|}
    \hline 
    Params. & Description & Expression/Values & References  \\
    \hline
    \hline
    $\kappa$ & Chem. diffusion rate & $10^{-6} - 10^{-5}$ $\text{cm}^2/\text{s}$ & \citealt{Stocker2012} \\
    $\beta$ & Chem. decay rate & $5-40$ $\text{s}^{-1}$ &  \citealt{Kawaguchi2011}\\
    $\lambda$& Tumbling rate & $2-50$  $\text{s}^{-1}$ & \citealt{Visser2003}\\
    $v$ & Swimming speed & $1 \times 10^{-2}$ cm/s & \citealt{Harvey2015} \\
    \hline
\end{tabular}
\caption{Dimensional parameters with descriptions and experimental values, along with citations.}
\label{tab:ParametersD}
\end{table}
\begin{table}[t!]
\centering
\renewcommand\arraystretch{1.3}
\begin{tabular}{|c|c|c|c|}
\hline
    Parameter & Description & Important Details & Simulation Values \\
    \hline
    \hline
    $d_1$ & Chem. diffusion rate  & $\frac{\kappa \lambda_0}{v^2}$ & $0.1 - 5$ \\
    $d_2$ & Chem. decay rate & $\frac{\beta}{\lambda_0}$ &  $0.1 - 5$ \\
    $\delta$ & Run-and-Tumble parameter & -- &  $0 - 0.02$ \\
    $\ell$ & Length of domain & -- & $3-10$ \\
    $\gamma$ & Max. chemical strength & -- & $0.01-0.05$\\
    $c_0$ & Chem. sensitivity threshold & $f_2$ and $f_3$ only & $0.005-0.05$\\
    $T_0$ & Behavior switch length & -- & $0.0001 - 0.05$\\
    $N_p$ & Number of plankton & 2D system only & $10^5 - 10^6$ \\
    \hline
\end{tabular}
\caption{Non-dimensional parameters used in 1D and 2D simulations along with descriptions, expressions of values in terms of dimensional parameters, and important details, along with the numerical values used in our simulations.}
\label{tab:ParametersND}
\end{table}

\section{Analysis of the 1D auto-chemotactic system}\label{sec:Analysis1D}
\subsection{Stability Analysis}\label{sec:Analytic1D}
To understand the dynamics of the 1D system, we analyze the stability of the constant steady states of $\eqref{eq:nd1}-\eqref{eq:nd3}$. There is a constant steady state from $\eqref{eq:nd1}$ of $d_2 \overline{c} = f(\overline{c})\overline{\rho}$, where $\overline{\rho} = 2\overline{\psi}$ and $\overline{\psi}$ is the constant solution for both $\psi^+$ and $\psi^-$. We linearize around this solution by letting 
\begin{equation}\label{eq:Perturbs1D}
    c = \overline{c} + \epsilon b, \quad \psi^+ = \overline{\psi} + \epsilon a^+, \quad \psi^- = \overline{\psi} + \epsilon a^-,
\end{equation} 
where $\epsilon \ll 1$ and $a^+, a^-,$ and $b$ are functions in terms of time and space. Substituting \eqref{eq:Perturbs1D} into $\eqref{eq:nd1}-\eqref{eq:nd3}$, at leading order of $\epsilon$ we obtain
\begin{dgroup}
\begin{dmath}\label{eq:1storder1}
b_t = d_1 b_{xx} + d_3 b + f(\overline{c})(a^+ + a^-)
\end{dmath}
\begin{dmath}\label{eq:1storder2}
a_t^+ = -a_x^+ + \frac{1}{2}(a^- - a^+) + \frac{\overline{\psi} b_x}{\delta} 
\end{dmath}
\begin{dmath}\label{eq:1storder3}
a_t^- = a_x^- - \frac{1}{2}(a^- - a^+) - \frac{\overline{\psi} b_x}{\delta}, 
\end{dmath}
\end{dgroup}
where $d_3 = \overline{\rho} f'(\overline{c}) - d_2$. As a note, we utilized the Taylor expansion $\frac{x}{\sqrt{x^2 + \delta^2}} = \frac{x}{\delta} - \frac{x^3}{2 \delta^2} + \mathcal{O}(x^5)$. We now consider the Fourier transform of the \eqref{eq:1storder1}-\eqref{eq:1storder3}, and let $A^+$, $A^-$, and $B$ denote the transforms of $a^+$, $a^-$ and $b$, respectively. Letting $k$ denote the wave number in the Fourier transform, we obtain
\begin{dgroup}
\begin{dmath}\label{eq:FT1}
B' = -k^2 d_1 B + d_3 B + f(\overline{c})(A^+ + A^-)
\end{dmath}
\begin{dmath}\label{eq:FT2}
(A^+)' = -ik A^+ + \frac{1}{2}(A^- - A^+) + \frac{ik \overline{\psi}}{\delta} B
\end{dmath}
\begin{dmath}\label{eq:FT3}
(A^-)' = ik A^+ - \frac{1}{2}(A^- - A^+) - \frac{ik \overline{\psi}}{\delta} B,
\end{dmath}
\end{dgroup}
which can be rewritten into the matrix form 
\begin{equation}\label{eq:StableMatrix}
    \begin{bmatrix}
    B \\
    A^+ \\
    A^- 
    \end{bmatrix}'
    = 
    \begin{bmatrix}
       d_3 - d_1 k^2 & f(\overline{c}) & f(\overline{c})  \\
    \frac{ik \overline{\psi}}{\delta}  & -ik - \frac{1}{2}  & \frac{1}{2} \\
     - \frac{ik \overline{\psi}}{\delta} & \frac{1}{2} & ik - \frac{1}{2}
    \end{bmatrix}
    \begin{bmatrix}
    B \\
    A^+ \\
    A^-
    \end{bmatrix}.
\end{equation}
We now look for eigenvalues of the $3 \times 3$ matrix in \eqref{eq:StableMatrix}. Denoting $\lambda$ as an eigenvalue of the above matrix, we can write the characteristic equation as

\begin{dmath}\label{eq:charact}
     \lambda^3 +\left[ d_1 k^2 - d_3 + 1\right] \lambda^2  + \left[ (d_1 + 1)k^2 - d_3 \right] \lambda  + k^2 \left( d_1 k^2 - d_3 - \frac{\overline{c} d_2}{\delta}\right) = 0.
\end{dmath}

\begin{figure}[t!]
    \centering
    \subfloat[Varying $d_1$ with $d_2 = 1$]{
        \includegraphics[width=0.48\linewidth]{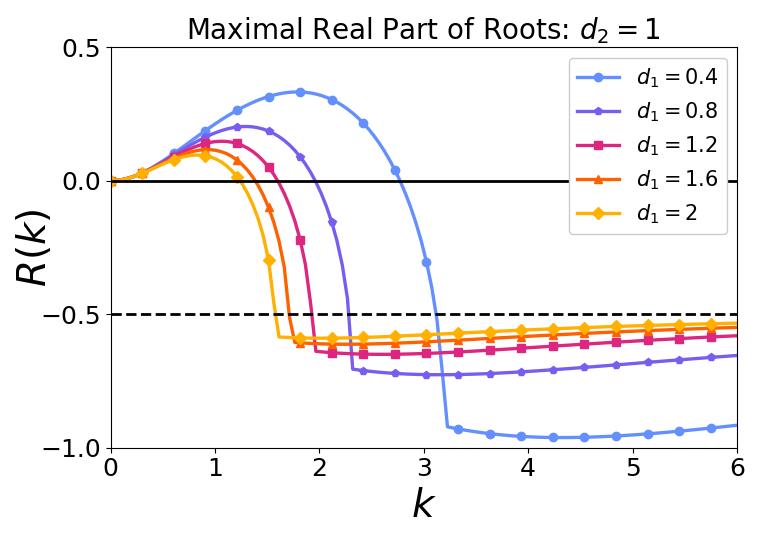}\label{fig:MaxReald1}}
        \hspace{0.01\linewidth}
    \subfloat[Varying $d_2$ with $d_1 = 1$]{
        \includegraphics[width=0.48\linewidth]{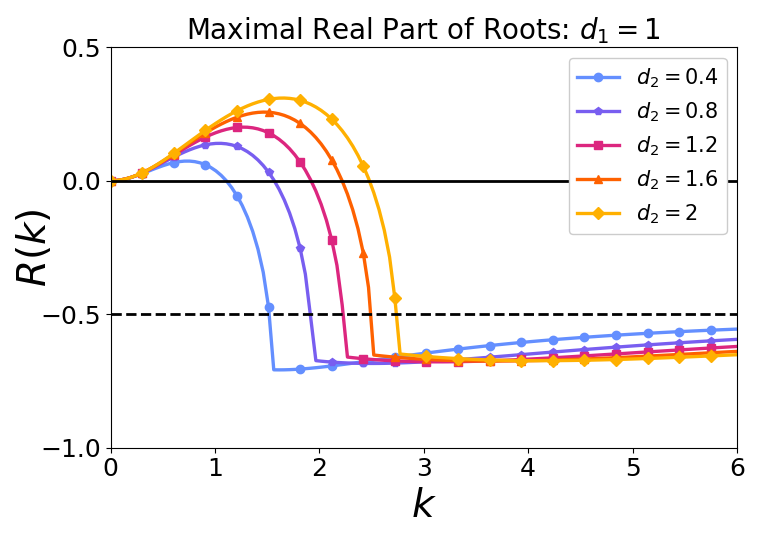}\label{fig:MaxReald2}}\\
    \subfloat[Stability regions with $\delta = 0.01$]{
        \includegraphics[width=0.48\linewidth]{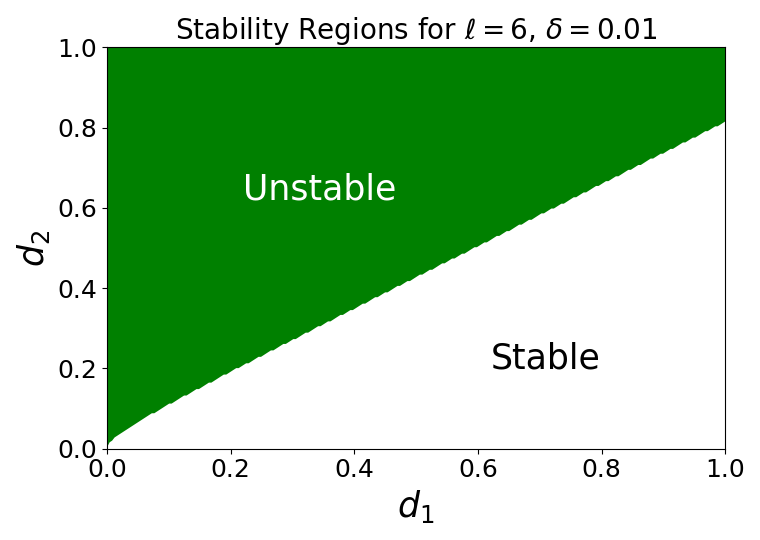}\label{fig:1DStabled1}} 
        \hspace{0.01\linewidth}
    \subfloat[Stability regions with $\delta = 0.012$]{
        \includegraphics[width=0.48\linewidth]{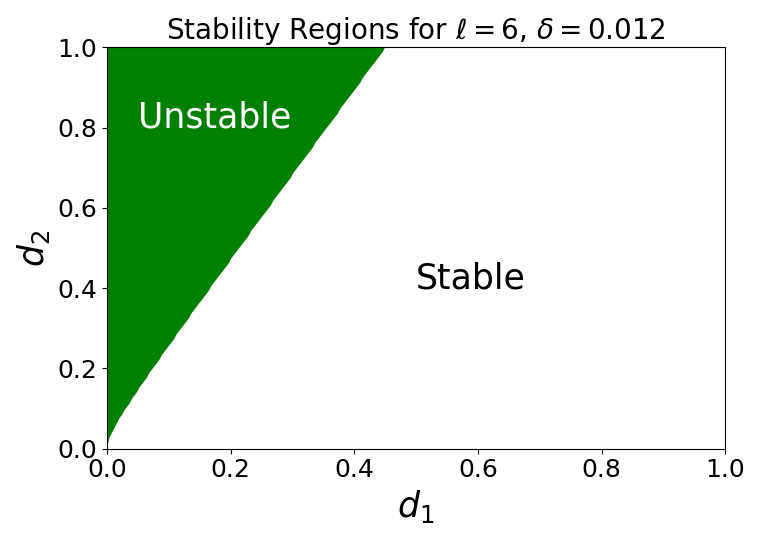}\label{fig:1DStabled2}}
    \caption{Plots \textbf{a} and \textbf{b} display $R(k)$ with varied $d_1$ and $d_2$ parameters. Plots \textbf{c} and \textbf{d} show phase space stability for a constant solution of \eqref{eq:nd11}-\eqref{eq:nd12} with varied $\delta$. Parameters used can be found in Appendix \ref{sec:1DParams}.} 
    \label{fig:RealRoot}
\end{figure}

We seek to find combinations of parameters in which the constant solution is unstable. Take the roots of \eqref{eq:charact} as the set $\{ \lambda_1(k), \lambda_2(k), \lambda_3(k)\}$ and define $R(k) := \max \{\text{Re}\left(\lambda_i\right)\}_{i =1}^3$. Let $k_u$ be the most unstable wave number, i.e.
\begin{equation}\label{eq:MostUnst}
k_u = \argmax\limits_{k > 0} R(k).
\end{equation}
If $k_u > \frac{2\pi}{\ell}$, the set of parameters will make $(\overline{c}, \overline{\rho})$ unstable in the domain length $\ell$. Given a fixed set of parameters $d_1, d_2, \overline{c}, \delta$ and deposition function $f(c)$, we can vary $k$ to numerically find the most unstable wave number. Examples of the function $R(k)$ are shown in Fig. \ref{fig:RealRoot}. In Fig. \ref{fig:MaxReald1}, we keep $d_2$ constant but vary the diffusion rate, $d_1$. As $d_1$ increases, $k_u$ decreases. Alternatively, in Fig. \ref{fig:MaxReald2}, $d_1$ is kept constant but the chemical decay rate, $d_2$, is varied. As $d_2$ increases, $k_u$ increases.

To see this relationship between $d_1$ and $d_2$ more clearly, we plot the stability regions for a domain length $\ell$ in the $d_1$-$d_2$ plane, shown in Figs. \ref{fig:1DStabled1} and \ref{fig:1DStabled2}. The green, shaded regions are examples of parameter regimes where the system would be unstable when $\ell = 6$ for either $\delta = 0.01$ or $\delta = 0.012$. As $\delta$ increases incrementally, the stable regime becomes much larger. As well, we note that the curve separating the unstable and stable regimes becomes linear as $d_1$ and $d_2$ both increase. 

We now remark on some properties of \eqref{eq:charact} and $R(k)$, the proofs of which can be found in Appendices \ref{sec:Proof1} and \ref{sec:Proof2}, respectively. These two propositions give stability conditions for large wave numbers and when the constant solution $(\overline{c}, \overline{\rho})$ may admit an unstable solution. 
\begin{proposition}\label{prop1}
If $d_1 > 0$, $\lim\limits_{k \to \infty} R(k) = -\frac{1}{2}$.
\end{proposition}
\begin{proposition}\label{prop2}
Define $d_3 = \overline{\rho} f'(\overline{c}) - d_2$. Then, if 
\begin{equation}\label{eq:ConditionProp2}
    -\frac{\overline{c} d_2}{\delta} < d_3 < \frac{1}{2d_1} \left[1 + \sqrt{1 + \frac{4 \overline{c} d_1 d_2}{\delta}} \right],
\end{equation}
there will be an $\ell$ such that the constant solution $(\overline{c}, \overline{\rho})$ to $\eqref{eq:nd1}-\eqref{eq:nd3}$ is unstable. 
\end{proposition}

\subsection{Linear Growth and Nonlinear Saturation}\label{sec:NumRes1D}
Having analyzed regimes which admit unstable or stable constant solutions, we now explore the long term dynamics of the model. We computed the nonlinear solutions shown in this section using numerical simulations of our one dimensional model represented by \eqref{eq:nd11}-\eqref{eq:nd12} with periodic boundary conditions. The specific details of our algorithm and the parameters used are in Appendix \ref{sec:1DNum} and Table \ref{tab:1DParams}, respectively.

\begin{figure}[t!]
    \centering
        \includegraphics[width=\linewidth]{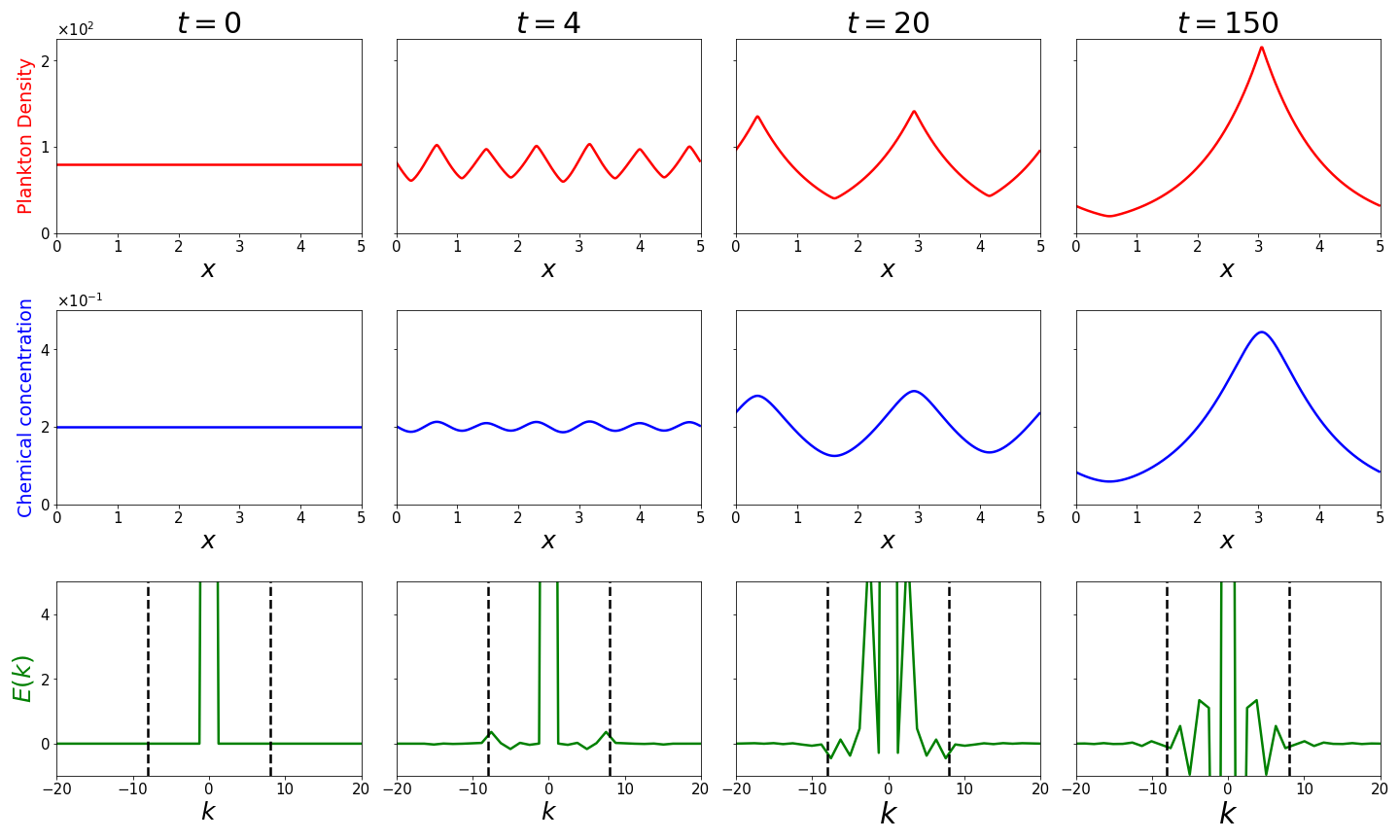}
    \caption{Time-series of plots for simulation of 1D autochemotactic system \eqref{eq:nd11}-\eqref{eq:nd12} at $t = 0, 4, 20, 150$. The top row displays the evolution of plankton density, $\rho$. The middle row display the evolution of chemical concentration, $c$. The bottom row shows the evolution of $E(k)$ (see \eqref{eq:FT1DEq}), where $k$ are the Fourier modes. The black vertical lines in the bottom row are at the most unstable wave number, $k_u \approx 8$. This simulation can be seen in its entirely in Online Resource 1 and parameters used can be found in Appendix \ref{sec:1DParams}.}
    \label{fig:SSNOnConst}
\end{figure}

In Figure \ref{fig:SSNOnConst} and Online Resource 1,  we observe the evolution of the system in a parameter regime where the constant stationary solution is unstable. The unstable distribution of aggregations rapidly saturates from nonlinear interactions and yields a slowly evolving coarsening pattern, eventually becoming a single aggregation. The top two panels display the plankton density, $\rho$, and chemical density, $c$, respectively. The lower panel displays $E(k)$, which shows the real part of the Fourier modes of $c$, i.e.
\begin{equation}\label{eq:FT1DEq}
    E(k;t_0) = \text{Re}\left\lbrace \mathcal{F}\left[c(x,t_0)\right](k)\right\rbrace,
\end{equation}
where $\mathcal{F}$ is the Fourier transform operator. To calculate this numerically, we utilize a fast Fourier transform. Aside from the translational mode of $k = 0$, the most unstable wave number, $k_u$, will grow the quickest initially. In this simulation, the most unstable wave number is $k_u \approx 8$, denoted by the black, dashed lines on the plots.
At $t = 4$, six small aggregations are formed and the mode at $k = k_u$ has been triggered, as it is the most pronounced $k$ at the outset. At $t = 20$, the aggregations have coarsened down to only 2 aggregations, and $k_u$ has diminished from view due to this clustering. The final panel at $t = 150$ shows that the two aggregations have merged, forming one large aggregation. We emphasize the timescale difference between the coarsening from six to two aggregations and from two to one aggregations.

\begin{figure}[t!]
    \centering
    \subfloat[Simulation using $f(c) = f_1(c)$ (animation in Online Resource 1)]{
        \includegraphics[width=0.99\linewidth]{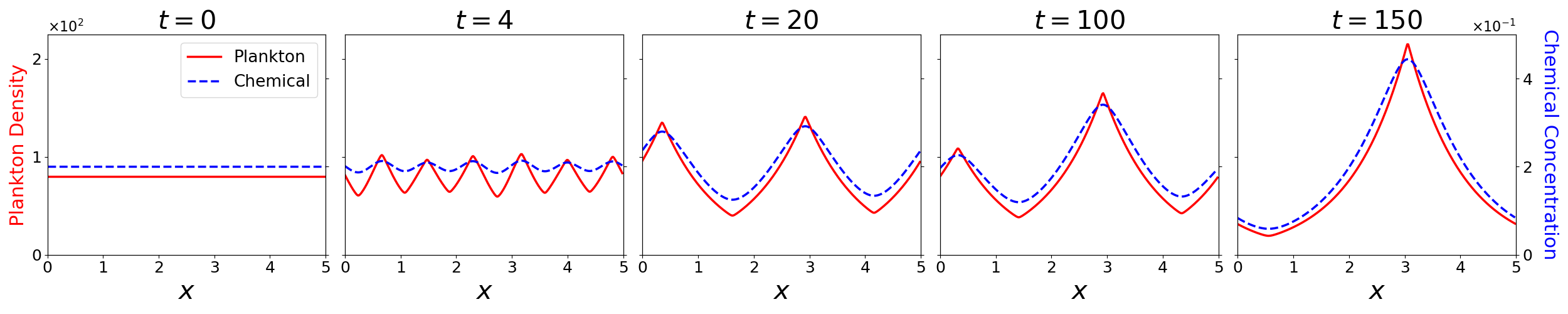}\label{fig:ConstantDep}}\\
        \subfloat[Simulation using $f(c) = f_2(c)$ (animation in Online Resource 2)]{
        \includegraphics[width=0.99\linewidth]{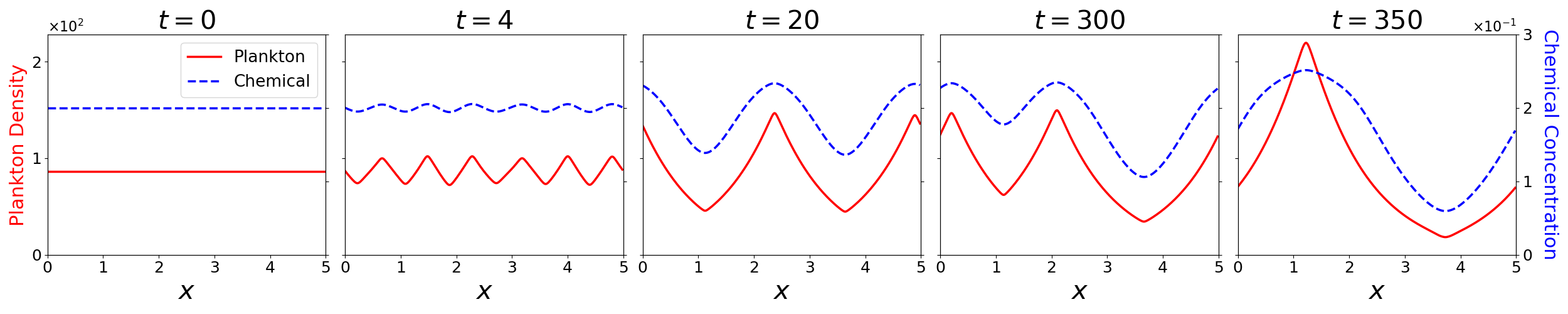}\label{fig:SwitchDep}}\\
    \subfloat[Simulation using $f(c) = f_3(c)$ (animation in Online Resource 3)]{
        \includegraphics[width=0.99\linewidth]{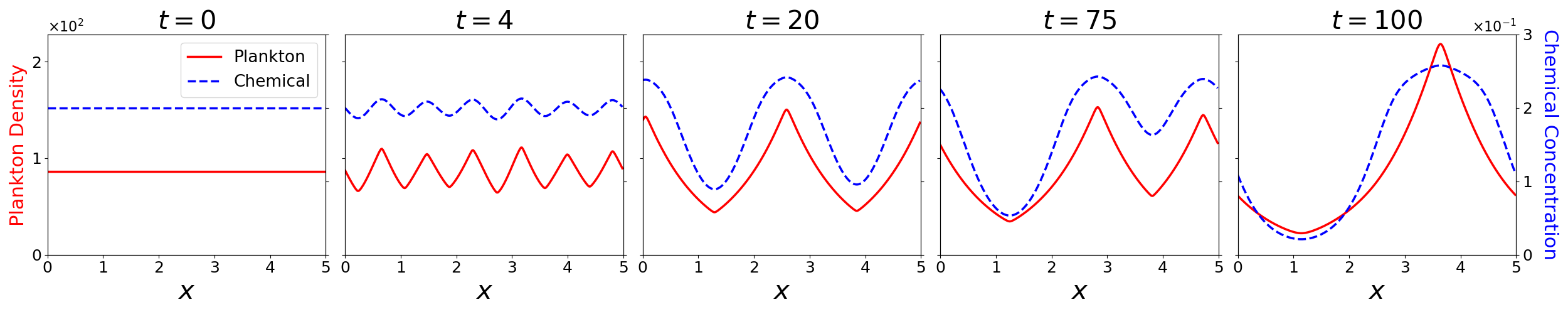}\label{fig:LinSwitchDep}} \\
            \subfloat[Total chemical over time, $C(t)$]{
        \includegraphics[width=0.5\linewidth]{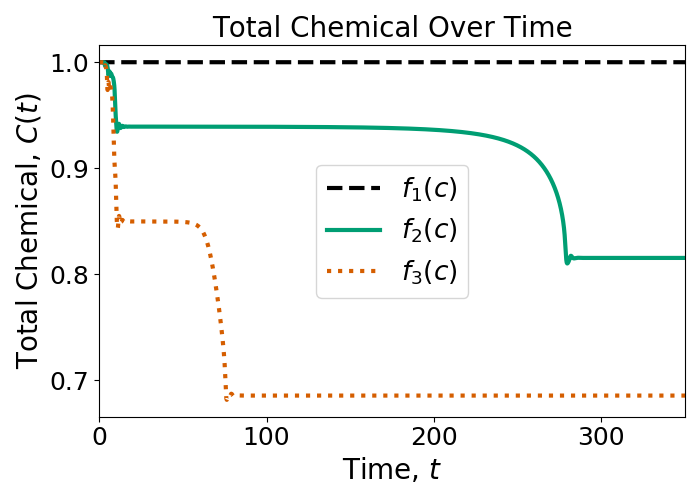}\label{fig:TotChem1d}} 

    \caption{Time-series of plots for simulation of 1D autochemotactic system \eqref{eq:nd11}-\eqref{eq:nd12} for various deposition functions. All simulation parameters, which can be found in Appendix \ref{sec:1DParams}, result in an unstable constant solution with $k_u = 8.169$. Note that the scale for the chemical concentration is different in \textbf{a}. The total chemical over time for all three simulations is in \textbf{d}.}
    \label{fig:EvolutionDeposition}
\end{figure}

\begin{figure}
    \centering
    \subfloat[Plankton density with varied $\delta$]{
        \includegraphics[width=0.45\linewidth]{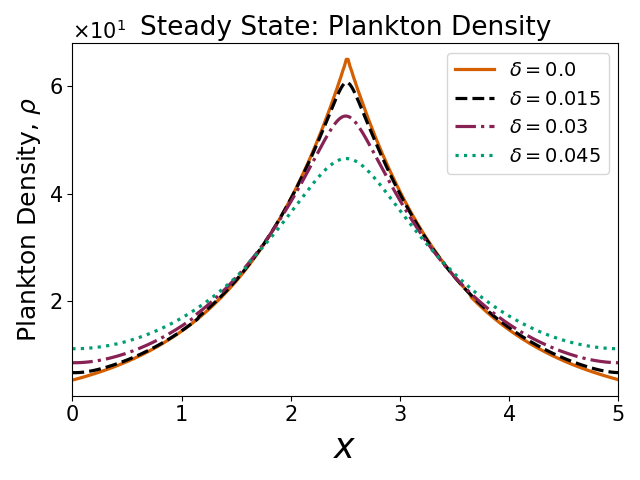}\label{fig:DeltaPlank}} \hspace{0.02\linewidth}
        \subfloat[Chemical concentration with varied $\delta$]{
        \includegraphics[width=0.45\linewidth]{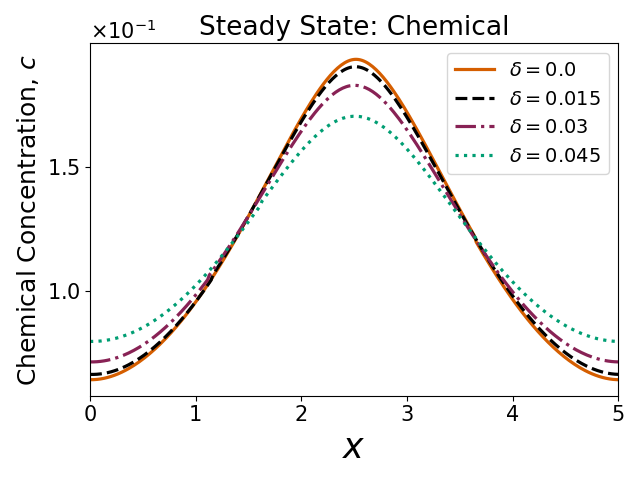}\label{fig:DeltaChem}}
    \caption{The non-constant steady states for the chemical concentration and plankton density for varied $\delta = 0, 0.015, 0.03, 0.045$. For comparison purposes, all steady states have been shifted to the center of the domain. Parameters can be found in Appendix \ref{sec:1DParams}.} 
    \label{fig:DeltaProg}
\end{figure}

In Figs. \ref{fig:EvolutionDeposition} and \ref{fig:DeltaProg}, we explore two important facets of our model: the deposition function and $\delta$. Recall our three deposition functions defined in \eqref{eq:Switch} and \eqref{eq:LinearSwitch}. For this comparison, we choose parameter values for $f_1(c)$, $f_2(c)$, and $f_3(c)$ such that all systems are unstable and have the same $k_u$.

Some key differences emerge from these simulations in Figs. \ref{fig:ConstantDep}-\ref{fig:LinSwitchDep}. Even though all simulations have a similar number of aggregations from $t = 4$ to $t = 20$, ending up with two aggregations at $t = 20$, the time needed to fully merge the final two aggregations varies. The linear switch deposition function, $f_3$, admits a singular aggregation at $t = 100$, but it takes until $t = 350$ for the switch deposition function, $f_2$, to merge. As well, we see a clear difference in the profile for $c$ achieved at the non-constant steady state. The $f_1$ function admits a profile for $c$ and $\rho$ that are relatively similar, but $c$ and $\rho$ look vastly different for the $f_2$ and $f_3$ cases. This can be attributed to the threshold parameter, $c_0$, hampering the plankton's ability to deposit chemical once $c$ becomes too large. 

Other insights can be gleaned from the total chemical concentration profile $C(t)$ (Fig. \ref{fig:TotChem1d}). The total chemical for the $f_1$ simulation remains at its equilibrium value. Between merging events for the $f_2$ and $f_3$ simulations, there is steady chemical total; however, several time steps prior to a merging event, there is a drastic decrease in the amount of chemical due to the restructuring of $\rho$. We show that $C(t)$ stays constant for $f(c) = f_1(c)$ in a general two-dimensional case in Proposition \ref{prop:Prop3} in Section \ref{sec:Analysis2D}.

This model is also highly sensitive to $\delta$, which expresses the ability of the organism to admit strong differentiated responses to small gradients. In Fig. \ref{fig:DeltaProg}, the $\delta$ parameter is decreased to show the differences in the final steady state solution of the model. As $\delta$ decreases, the peak of the plankton density and chemical concentration become more sharp. At $\delta = 0$, note that a steady state solution of \eqref{eq:nd11}-\eqref{eq:nd12} must satisfy 
\begin{equation}\label{eq:NonConstSS}
    \rho_{xx} - \left[ \text{sgn}(c_x) \rho \right]_x = 0,
\end{equation}
with periodic boundary conditions. Due to the translational invariance of the solution to \eqref{eq:NonConstSS}, we can center the aggregation in the middle of our domain $[0,\ell]$. Using Fig. \ref{fig:DeltaPlank}, we recognize that $\text{sgn}(c_x) > 0$ where $0  \leq x \leq \ell/2$ and $\text{sgn}(c_x) < 0$ where $\ell/2 \leq x \leq \ell$. Thus, we obtain a steady state solution of stitched exponential functions:
\begin{equation}\label{eq:SSf}
 \rho(x) = \begin{cases} A_0 e^{x}, & x \in [0,\ell/2] \\
A_0 e^{-(x - \ell)}, & x \in [\ell/2,\ell] 
\end{cases}, \quad A_0 = \frac{\overline{\rho} \ell}{2\left(e^{\ell/2} - 1\right)}.
\end{equation}
The stability of this steady state is still unknown, along with the other non-constant steady state solutions with $\delta > 0$. However, through numerical experimentation, we conjecture that this steady state is unconditionally stable.

\section{Analysis of the 2D auto-chemotactic system}\label{sec:Analysis2D}
\subsection{Fourier stability analysis}\label{sec:Analytic2D}
We now look for steady state solutions to the two-dimensional system \eqref{eq:2nond1}-\eqref{eq:2nond3}. We find a trivial constant steady state of $d_2 \overline{c} = f(\overline{c}) \overline{\rho}$ and $\overline{\rho} = 2 \pi \overline{\psi}$, similar to the one-dimensional analysis. We perturb this constant solution as
\begin{equation}\label{eq:perturb}
    \psi = \overline{\psi} + \epsilon \phi, \quad c = \overline{c} + \epsilon \mu , \quad \rho = \overline{\rho} + \epsilon \int_0^{2\pi} \phi := \overline{\rho} + \epsilon \tau, 
\end{equation}
where $\phi, \mu, \tau$ are functions of time and space and $\epsilon \ll 1$. Substituting the perturbations from \eqref{eq:perturb} to \eqref{eq:2nond1}, we obtain the evolution equation for $\mu$ at first order as
\begin{equation}\label{eq:chempert}
    \mu_t = d_1\Delta \mu + d_3 \mu + f(\overline{c}) \tau,
\end{equation}
where $d_3 = f'(\overline{c}) \overline{\rho}-d_2$. By substituting the perturbations into \eqref{eq:2nond3} and utilizing the same Taylor expansion as in the one-dimensional case for the terms with $\delta$, we obtain
\begin{dgroup*}
\begin{dmath}\label{eq:MidStepPert}
\epsilon \phi_t = - \epsilon \vec{e}_{\theta} \cdot \nabla \phi - \frac{1}{2}\left( 1 - \epsilon \frac{\nabla \mu \cdot \vec{e}_{\theta}}{\delta} \right)(\overline{\psi} + \epsilon \phi) + \frac{1}{4 \pi} \int_0^{2\pi} \left[ 1 - \epsilon \frac{\nabla \mu \cdot \vec{e}_{\theta'}}{\delta} \right] (\overline{\psi} + \epsilon \phi(\theta') ) \, d\theta'.
\end{dmath}
\end{dgroup*}
At first order, we simplify \eqref{eq:MidStepPert} to become
\begin{equation} \label{eq:plankpert}
\phi_t = - \vec{e}_{\theta} \cdot \nabla \phi - \frac{1}{2}\phi + \frac{\overline{\psi}}{2\delta} \nabla \mu \cdot \vec{e}_{\theta} + \frac{1}{4\pi} \tau.
\end{equation}
Taking the two-dimensional Fourier transform in both $x$ and $y$ of \eqref{eq:chempert} and \eqref{eq:plankpert} while denoting $\vec{k} = \langle k_1, k_2 \rangle$ as the wave vector, we get
\begin{dgroup}
\begin{dmath}\label{eq:HFt}
H_t  = i (\vec{k} \cdot \vec{e}_{\theta}) H  - \frac{1}{2} H - \frac{\overline{\psi} i }{2 \delta}(\vec{k} \cdot \vec{e}_{\theta}) M + \frac{1}{4\pi} T
\end{dmath}
\begin{dmath}\label{eq:MFt}
M_t = (d_3 - d_1 |\vec{k}|^2) M + f(\overline{c}) T,
\end{dmath}
\end{dgroup}
where $H$, $M$, and $T$ are the two-dimensional Fourier transforms of $\phi$, $\mu$, and $\tau$, respectively. Continuing, we write $H$ in terms of its Fourier coefficients with respect to $\theta$, i.e.
\begin{equation}\label{eq:HCoeff}
     H(\vec{k}, t, \theta) = \sum\limits_{n \in \mathbb{Z}} h_n(\vec{k}, t) e^{in \theta},
\end{equation}
where $h_n \in \mathbb{C}$. From this representation, $T$ can be written as
\begin{equation}\label{eq:TCoeff}
    T(\vec{k}, t) = \int_0^{2\pi} H(\vec{k}, t, \theta') \, d \theta' = 2 \pi h_0(\vec{k},t).
\end{equation}
We define the complexification of the wave number $\vec{k}$ as $\omega := k_1 + i k_2$, and denote $\overline{\omega}$ as the complex conjugate of $\omega$. We rewrite $\sin(\theta)$ and $\cos(\theta)$ in exponential forms and substitute \eqref{eq:HCoeff} and \eqref{eq:TCoeff} into \eqref{eq:HFt}-\eqref{eq:MFt}. By matching powers of $e^{i n \theta}$, we obtain the following infinite linear ODE system, where primes denote derivative with respect to $t$: 
\begin{dgroup}
\begin{dmath}
    h_0'  = \frac{i}{2}\left( \overline{\omega} h_{-1} + \omega h_1 \right) 
\end{dmath}
\begin{dmath}
    h_1'  = \frac{i}{2} \left( \overline{\omega} h_{0} + \omega h_2 \right) - \frac{1}{2} h_1 - \frac{\overline{\psi} i \overline{\omega}}{4 \delta} M
    \end{dmath}
    \begin{dmath}
     h_{-1}'  = \frac{i}{2} \left( \overline{\omega} h_{-2} + \omega h_0 \right) - \frac{1}{2} h_{-1} - \frac{\overline{\psi} i \omega}{4 \delta} M
     \end{dmath}
     \begin{dmath}
    h_{j}' =  \frac{i}{2} \left( \overline{\omega} h_{j-1} + \omega h_{j + 1} \right) - \frac{1}{2} h_j
    \end{dmath}
    \begin{dmath}
    M' = (d_3 - d_1 |\vec{k}|^2)M + 2 \pi f(\overline{c}) h_0,
    \end{dmath}
\end{dgroup}
where $j \neq -1, 0, 1$. This system can be rewritten as matrix equation with infinite entries by $\vec{y}' = \vec{A} \vec{y}$, where
\begin{equation}\label{eq:fourierMat}
\vec{y} = \begin{bmatrix}
\vdots \\
h_{-2}\\
h_{-1}\\
h_0\\
M\\
h_1 \\
h_2\\
\vdots 
\end{bmatrix}, \hspace{2mm} 
\vec{A} = 
\begin{bmatrix}
\ddots & \ddots  & \ddots & & & & & &  \\
 & \frac{i}{2} \overline{\omega} & -\frac{1}{2} & \frac{i}{2} \omega & & & & & \\
 & &  \frac{i}{2} \overline{\omega} & -\frac{1}{2} &  \frac{i}{2} \omega & -\frac{\overline{\psi} i \omega}{4 \delta} & & & & \\
 &  & &  \frac{i}{2} \overline{\omega} & 0 & 0 & \frac{i}{2}\omega & & & \\
 & &  & & 2 \pi f(\overline{c}) & d_3 - d_1 |\vec{k}|^2 & 0 & & &  \\
  & &  & &  \frac{i}{2} \overline{\omega} &- \frac{\overline{\psi} i \overline{\omega}}{4 \delta} & -\frac{1}{2} & \frac{i}{2}\omega & &  \\
   & &  & &  &0 & \frac{i}{2} \overline{\omega} & -\frac{1}{2} & \frac{i}{2}\omega &  \\
   & &  & &  &&  & \ddots  & \ddots  & \ddots  \\
\end{bmatrix}
\end{equation}
We now look for the eigenvalues of $\vec{A}$ which have the largest real part. Since we cannot compute every eigenvalue of $\vec{A}$ algebraically, we will compute them numerically by truncating the matrix at an appropriate threshold. Denote $\vec{A}_N$ and $\vec{y}_N$ as the truncated matrix $\vec{A}$ and vector $\vec{y}$ in \eqref{eq:fourierMat} consisting of the equations with the functions $\{ h_{-N}, \cdots h_{-1}, h_0, M, h_1, \cdots, h_N\}$, where $N \in \mathbb{Z}^+$ is sufficiently large (see next paragraph). We now calculate the eigenvalues of the system $\vec{y}_N' = \vec{A}_N \vec{y}_N$. Denote the set of these eigenvalues as $\left\lbrace \lambda_{i}(|\vec{k}|) \right\rbrace_{i = 1}^{2N + 2}$, and define $R_N = R_N(|\vec{k}|) := \max\lbrace \text{Re}\left(\lambda_i\right) \rbrace_{i = 1}^{2N + 2}$. Then, $K_u$ will be denoted as the norm of the most unstable wave number, i.e.
\begin{equation}\label{eq:MostUnst2d}
K_u := \left| \argmax\limits_{\vec{k}} R_N(|\vec{k}|) \right|.
\end{equation}
If $K_u \geq \frac{2\pi}{\ell}$, the set of parameters will make the constant solution $(\overline{c}, \overline{\rho})$ unstable in the domain $[0, \ell] \times [0, \ell]$. Given a fixed set of parameters $d_1, d_2, \overline{c}, \delta$ and deposition function $f(c)$, we vary $\vec{k}$ to find $K_u$.

\begin{figure}[p!]
    \centering
    \subfloat[Varying $d_1$ with $d_2 = 1$]{
        \includegraphics[width=.48\linewidth]{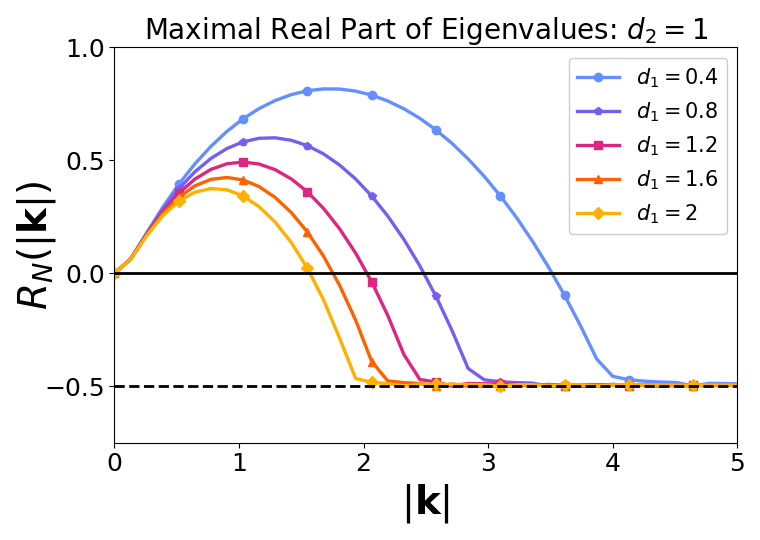}\label{fig:FirstUnstable2d}}
        \hspace{.01\linewidth}
    \subfloat[Varying $d_2$ with $d_1 = 1$]{
        \includegraphics[width=.48\linewidth]{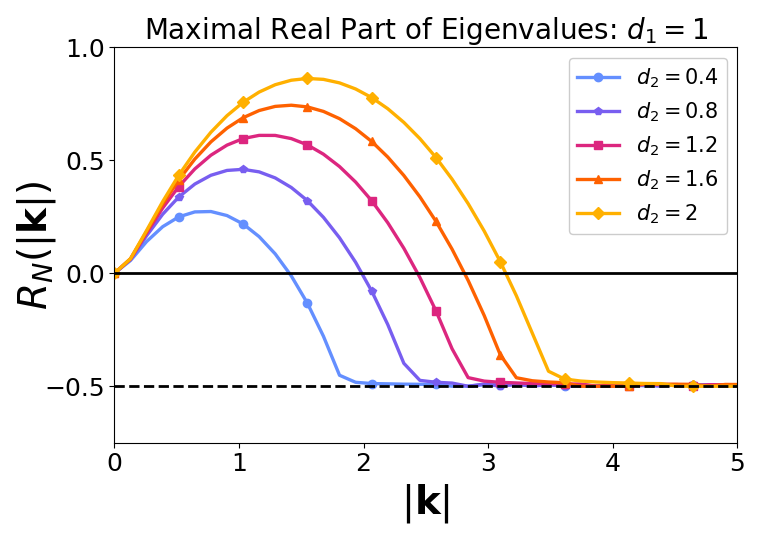}\label{fig:SecondUnstable2d}}\\
        \subfloat[Varying $N$ with $d_1 = d_2 = 1$]{
        \includegraphics[width=.48\linewidth]{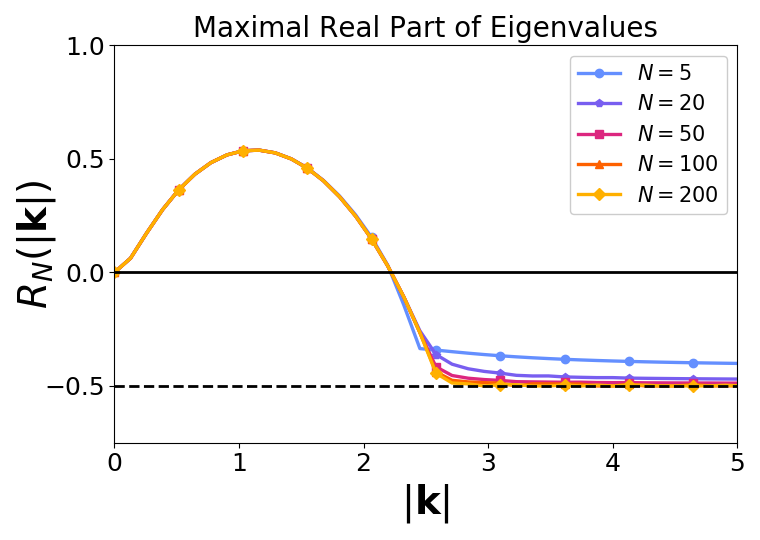}\label{fig:Ncheck}}
        \hspace{.01\linewidth}
    \subfloat[Varying $c_0$ with $d_1 = d_2 = 1$]{
        \includegraphics[width=.48\linewidth]{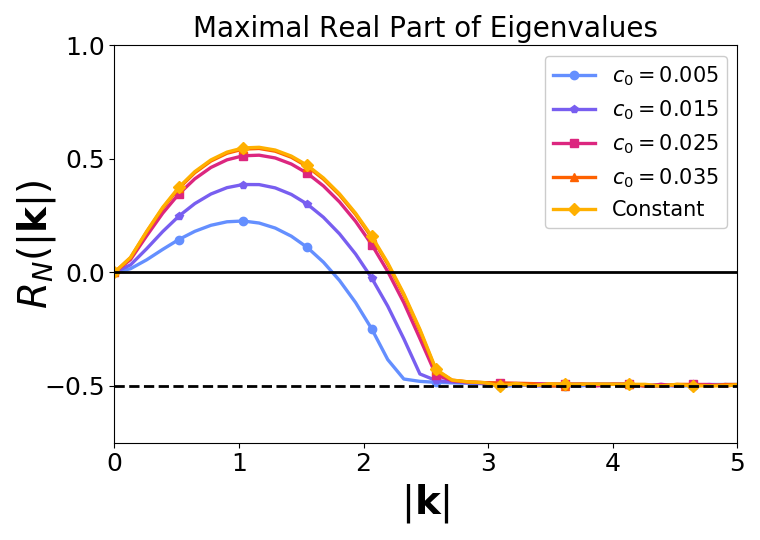}\label{fig:c0move}}\\
        \subfloat[Stability regions with $\delta = 0.001$]{
        \includegraphics[width=.48\linewidth]{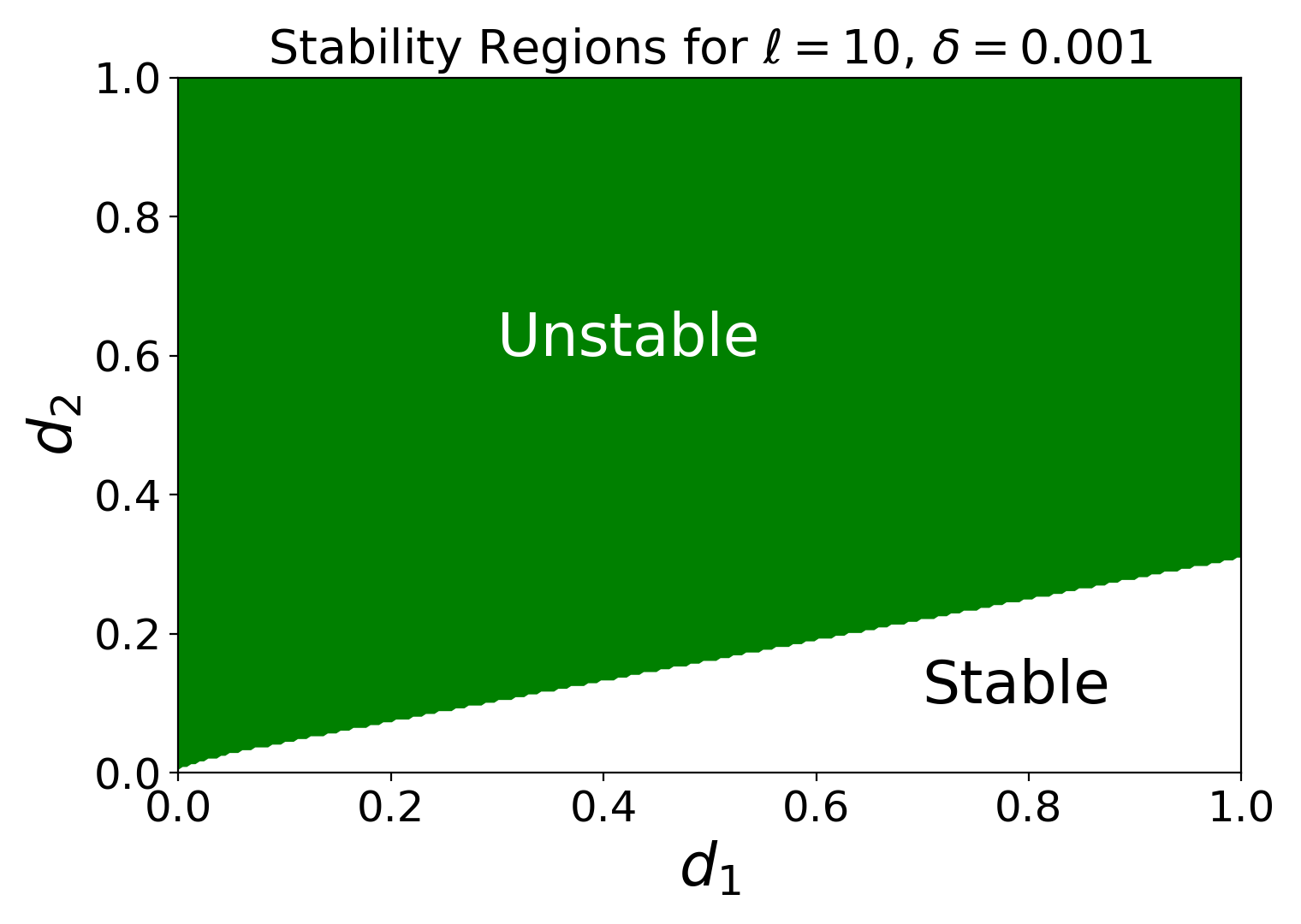}\label{fig:2DStability1delta}}
        \hspace{.01\linewidth}
    \subfloat[Stability regions with $\delta = 0.004$]{
        \includegraphics[width=.48\linewidth]{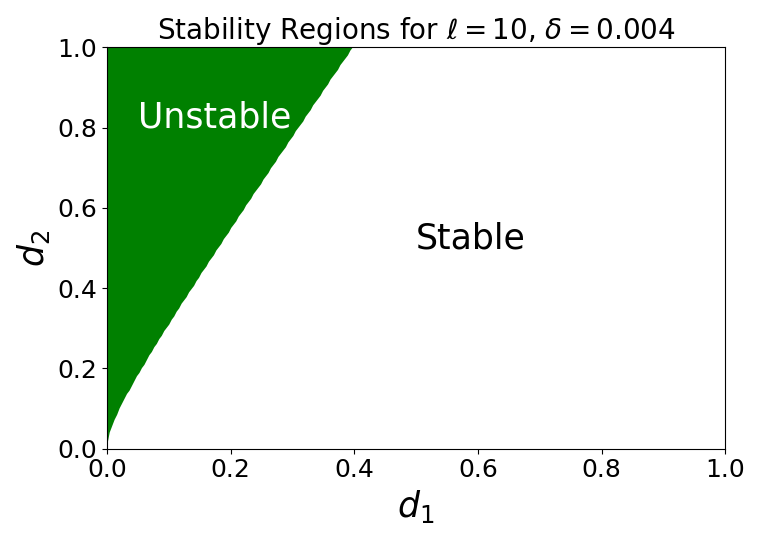}\label{fig:2DStability2delta}}
    \caption{Plots \textbf{a}-\textbf{d} display the function $R_N(|\vec{k}|)$, varying $d_1$, $d_2$, $N$ and $c_0$, respectively. For \textbf{a}, \textbf{b}, and \textbf{d}, we utilize $N = 100$. Plots \textbf{e} and \textbf{f} show phase space stability for a constant solution of \eqref{eq:2nond1}-\eqref{eq:2nond3} with varied $\delta$. Parameters used can be found in Appendix \ref{sec:2DParams}.}
    \label{fig:RealRoot2D}
\end{figure}

Examples of the function $R_N(|\vec{k}|)$ are shown in Figs. \ref{fig:FirstUnstable2d} and \ref{fig:SecondUnstable2d}, where the parameters $d_1$ and $d_2$ are varied. In a similar fashion to the one-dimensional case (Fig. \ref{fig:RealRoot}), decreasing $d_1$ or increasing $d_2$ causes $K_u$ to increase and $R_N(|\vec{k}|)$ approaches $-1/2$ as  $|\vec{k}|$ approaches infinity. In Fig. \ref{fig:Ncheck}, we found that $N = 100$ fully resolved $R_N(|\vec{k}|)$. This shows the rapid convergence of small $|\vec{k}|$ for all $N$ and as $N$ increases, the approximation of large $|\vec{k}|$ behavior converges towards an asymptote at $-1/2$. We also highlight how $c_0$, the threshold parameter for the deposition functions, changes the stability of the constant solution. In Fig. \ref{fig:c0move}, by increasing $c_0$ for $f_2$ in \eqref{eq:Switch}, the function $R_N(|\vec{k}|)$ for a constant deposition function with the same parameters is recovered.

Finally, Figs. \ref{fig:2DStability1delta} and \ref{fig:2DStability2delta} display the stability regions in the $d_1$-$d_2$ plane for two different $\delta$ values: 0.001 and 0.004. As $\delta$ increases, there are fewer steady states which will be unstable for a certain domain size $[0, \ell]\times [0,\ell]$. As with the one-dimensional case, the curve separating the stable and unstable regimes becomes linear as $d_1$ and $d_2$ both increase. 

\subsection{Linear Growth and Nonlinear Saturation}\label{sec:NumRes2D}
We now turn towards understanding the nonlinear saturation of unstable steady state solutions to our two-dimensional model. We analyze the Lagrangian system for $\rho$ and $c$ and show numerical results. Table \ref{tab:2DParams} and Fig. \ref{fig:Deposition} show the parameters and three deposition function we will utilize in our simulations, respectively.

To initialize our system, the density $\rho$ can be constructed as a summation of $N_p$ point sources, 
\begin{equation}
    \rho(\vec{x}, t) = \sum\limits_{i = 1}^{N_p} \delta(\vec{x} - \vec{x}_i(t)),
\end{equation}
where $\vec{x}_i$ is the position of each plankton and $\delta(\vec{x})$ is the Dirac delta function. Initially, all organisms are randomly distributed throughout the domain and given a random orientation $\theta_i$, selected from a uniformly random distribution $[0,2\pi)$. In each iteration, plankton tumble with a non-dimensional probability of 
\begin{equation}
    \mathbb{P}\left[ \text{Tumble in } (t, t + \Delta t) \right] = \frac{1}{2} \left( 1 - \frac{\vec{e}_{\theta_i}\cdot \nabla c}{\sqrt{\left( \vec{e}_{\theta_i} \cdot \nabla c\right)^2 + \delta^2}} \right)\Delta t.
    \label{eq:ProbND}
\end{equation}
If a plankton tumbles, a random direction is then chosen utilizing the turning kernel $T(\theta_i^{\prime}, \theta_i) = \frac{1}{2\pi}$. Each particle then updates their position by 
\begin{equation}
    \label{eq:UpdatePos}
    \vec{x}_i(t + \Delta t) = \vec{x}_i(t) + \Delta t \left\langle \cos(\theta_i), \sin(\theta_i) \right\rangle. 
\end{equation}

For autochemotaxis, we calculate $f(c)$ at these new positions using our deposition function of choice. To deposit the chemical numerically, we place a two-dimensional Gaussian centered at $\vec{x}_i$ instead of using a Dirac delta function. We then evolve equation \eqref{eq:2nond1} by utilizing a Crank-Nicolson time-stepping method. See Appendix \ref{sec:2DNum} for more specific details regarding the simulations.

\begin{figure}[ht!]
    \centering
        \includegraphics[width=.99\linewidth]{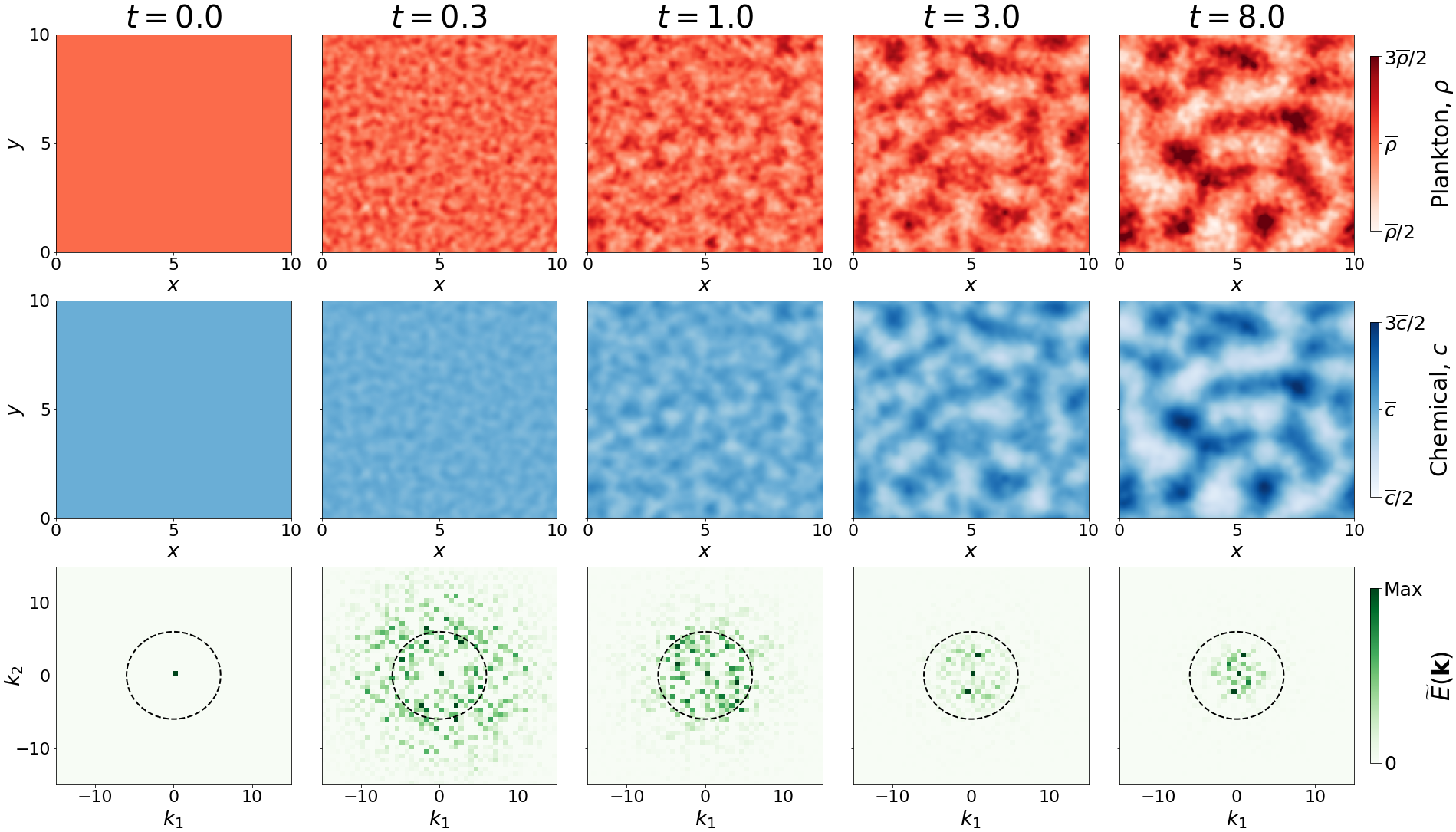}
    \caption{A time-series evolution of a 2D simulation described in Section \ref{sec:NumRes2D} at times $t = 0, 0.3,1,3,$ and $8$ for $\rho$, $c$, and $\widetilde{E}(\vec{k})$. Note that the scaling for $\widetilde{E}(k)$ changes every panel. The lightest shade is always 0, but the darkest shade is determined by the second largest value of $\widetilde{E}(\vec{k})$ at each $t$. The black, dotted circle in the bottom row has a radius $K_u$. For this simulation, $K_u = 6.096$. Parameters used can be found in Appendix \ref{sec:2DParams}. Animation of $\rho$ can be found on Online Resource 4.}
    \label{fig:2DEV}
\end{figure}

In Fig. \ref{fig:2DEV}, we display the evolution of a parameter regime which makes the constant solution $(\overline{c}, \overline{\rho})$ unstable and induce pattern formation. To visualize $\rho$ more effectively, we portray all $N_p$ plankton as small, two-dimensional Gaussians centered at $\vec{x}_i$ and plot the resulting density. The top row of Fig. \ref{fig:2DEV} shows this density evolution; the middle row displays the evolution of $c$, the chemical concentration; and the bottom row displays $\widetilde{E}(\vec{k})$, which calculates the real part of the Fourier modes of $c$, i.e.
\begin{equation}\label{eq:2DEtilde}
    \widetilde{E}(\vec{k};t_0) = \text{Re} \left\lbrace \mathcal{F}\left[ c(\vec{x}, t_0) \right](\vec{k})\right\rbrace,
\end{equation}
where $\mathcal{F}$ is the two-dimensional Fourier transform. To calculate this numerically, we utilize a two-dimensional fast Fourier transform. Other than the translational mode of $\vec{k} = \vec{0}$, a wave vector whose norm is equal to $K_u$ (see \eqref{eq:MostUnst2d}) should grow the quickest initially. For this simulation, the norm of the most unstable wave number is $K_u \approx 6$, denoted by the circle with radius $K_u$ in the lower panels.
At $t = 0.3$, the development of small aggregations occur and several modes grow quickly in the Fourier domain, the largest of which occurs at $K_u \approx 6$ as expected. As time progresses, the aggregations start to merge together and become more densely populated which can also be seen in $\widetilde{E}(\vec{k})$ as wave vectors with smaller norms emerge, reflecting pattern formation. Unsurprisingly, chemical concentrations are highest near the most densely populated areas of the plankton density.     
\begin{figure}[t!]
    \centering
        \includegraphics[width=.99\linewidth]{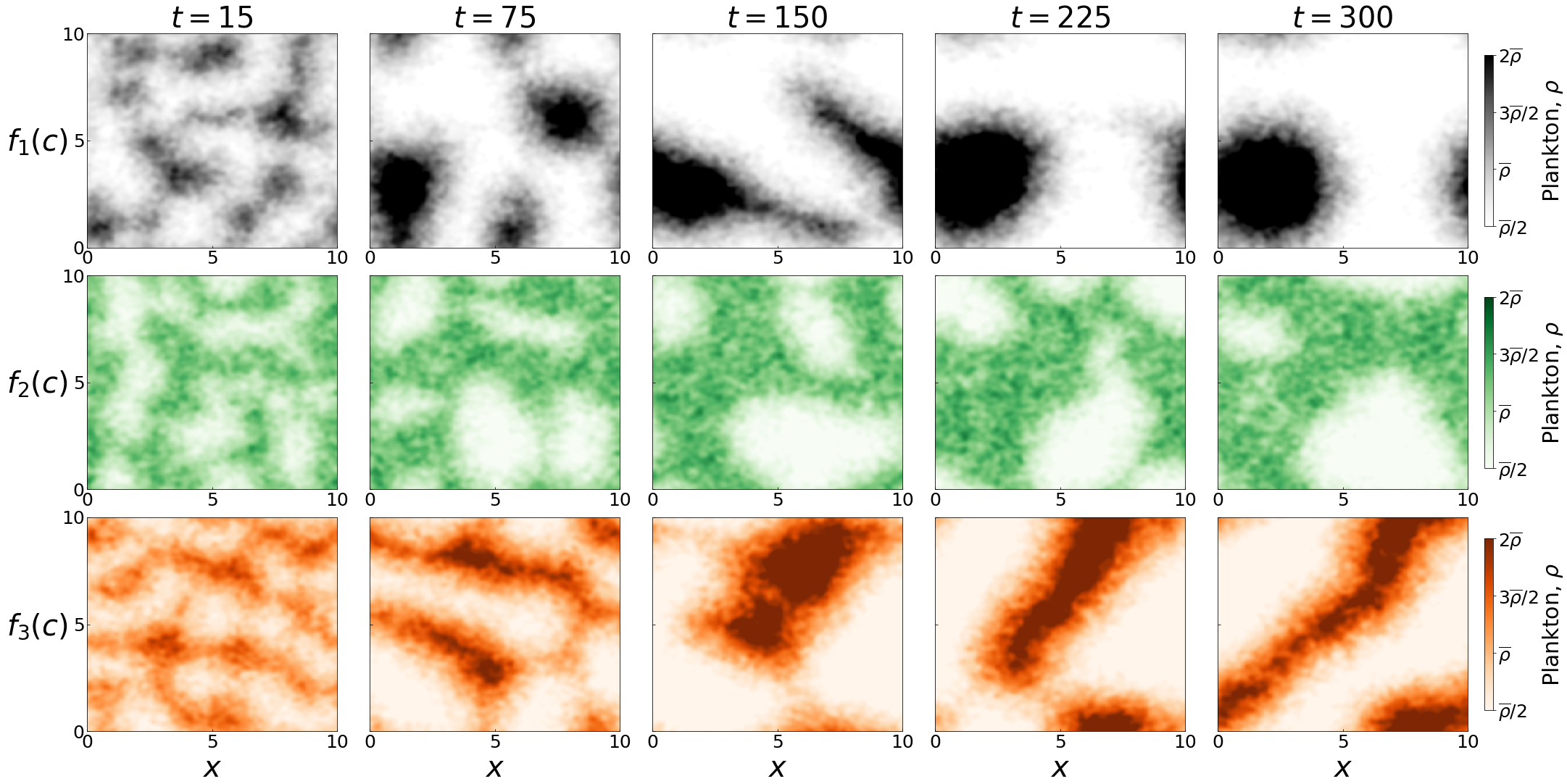}
    \caption{Time series evolution of the plankton density, $\rho$, described in Section \ref{sec:NumRes2D}, with varied deposition functions for times $t = 15, 75, 150, 225,$ and $300$. The top row is $f_1(c)$, the second row is $f_2(c)$ and the third row is $f_3(c)$, defined in \ref{eq:Switch} and \ref{eq:LinearSwitch}. All simulation parameters, which can be found on Table \ref{tab:2DParams} in Appendix \ref{sec:2DParams}, are constructed to have $K_u \approx 6$. Animations of these three simulations can be found on Online Resources 4, 5, and 6, respectively.}
    \label{fig:PlankDens}
\end{figure}

We now discuss the effects the various deposition functions have on the development and structure of aggregations. To equitably compare all three deposition functions, parameters were found such that $K_u \approx 6$ while keeping $d_1$, $d_2$, and $\delta$ the same (see Table \ref{tab:2DParams}). The evolution of the plankton densities can be seen in Fig. \ref{fig:PlankDens}. There are few differences in the initial development of the aggregations, as seen from the $t = 15$ plots. However, starting at $t = 75$, we qualitatively see marked differences. The constant deposition function, $f_1$, has more defined and densely populated aggregations, and at long times, it forms a singular dense and symmetric aggregation. The switch function, $f_2$,  The linear switch function, $f_3$, does have very populated areas throughout the domain and ends with a singular aggregation, but the aggregations are much less symmetric and have snake-like features. For all simulations, merging events occur more frequently early on rather than later, due to the aggregations' sizes.

\begin{figure}[t!]
    \centering
    \subfloat[Total chemical, $C(t)$]{
        \includegraphics[width=.48\linewidth]{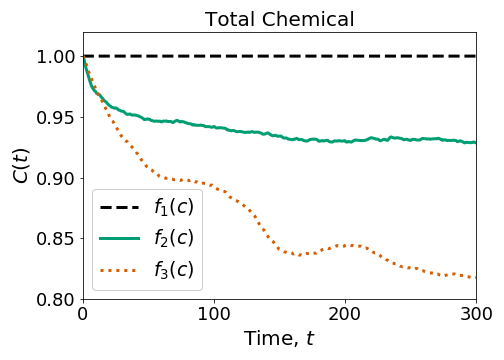}\label{fig:ChemicalTot}}
        \hspace{.01\linewidth}
    \subfloat[Difference in chemical, $\max(c(t))- \min(c(t))$]{
        \includegraphics[width=.48\linewidth]{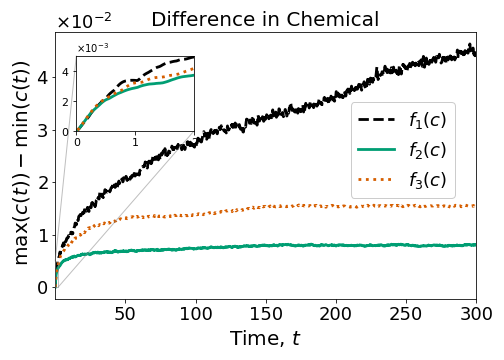}\label{fig:MaxMinChem}}\\
    \subfloat[Radial distribution function, $g(r)$]{
        \includegraphics[width=.99\linewidth]{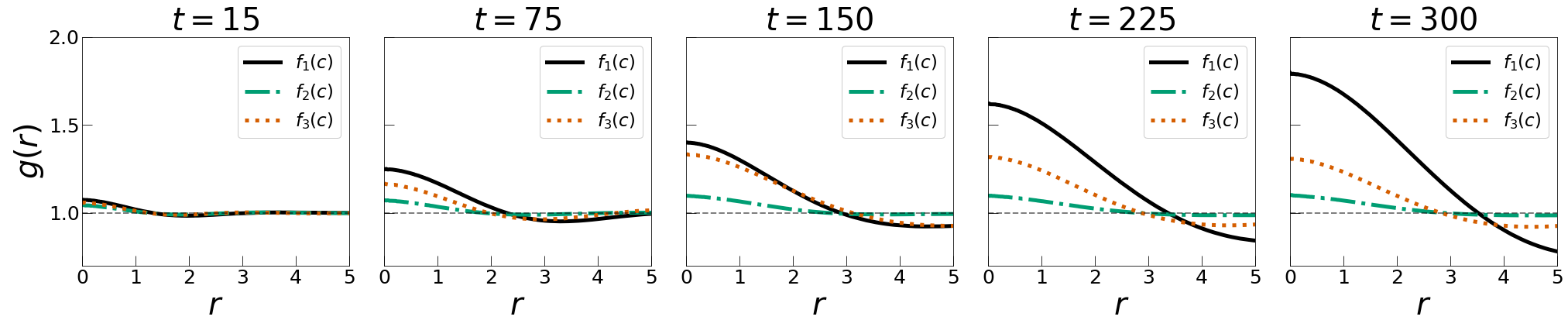}\label{fig:RDFs}}
    \caption{Plots \textbf{a}-\textbf{c} display various metrics to differentiate between simulations of the various deposition functions. For \textbf{a}, the total chemical, $C(t)$, is calculated \eqref{eq:TotChemeq}.  For \textbf{c}, the radial distribution function, $g(r)$, is defined in \eqref{eq:RDF}. All plots are averaged over six different simulations. Parameters for these simulations can be found in Appendix \ref{sec:2DParams}.}
    \label{fig:MetricsDeps}
\end{figure}

We utilize three metrics to elucidate the differences between these deposition functions, the first of which is the total chemical in the system. Denote $C(t)$ as the total chemical, i.e.
\begin{equation}\label{eq:TotChemeq}
    C(t) = \iint\limits_{D} c(\vec{x},t) \, d\vec{x},
\end{equation}
where $D = [0,\ell] \times [0,\ell].$ In Fig. \ref{fig:ChemicalTot}, the evolution of $C(t)$ is shown for the three simulations. All simulations were initialized with $c(\vec{x}, 0) = \overline{c} = 0.01$ and thus, $C(0) = \overline{c} \ell^2 = 1$. The long term behavior of all three simulations is decidedly varied, with $f_1$ staying constant, $f_2$ sharply declining before leveling out, and $f_3$ having periods of plateaus and sharp declines. Proposition \ref{prop:Prop3}, the proof of which is in Appendix \ref{sec:Proof3}, shows that the constant deposition function admits no change in chemical if the chemical is initialized at $\overline{c}$.
\begin{proposition}\label{prop:Prop3}
Denote the total chemical concentration as $C(t)$, defined in \eqref{eq:TotChemeq}. Let $d_2 \overline{c} = f(\overline{c}) \overline{\rho}$ be the constant steady state solution for \eqref{eq:2nond1}-\eqref{eq:2nond3} on $D = [0,\ell] \times [0,\ell]$. Then, if $f(c) = f_1(c)$ and $C(0) = C_0$ where $C_0 \geq 0$, then
\begin{equation}\label{eq:Prop3Eq}
    \lim\limits_{t \to \infty} C(t) = \overline{c} \ell^2.
\end{equation}
\end{proposition}

The second metric we utilize is the difference between the maximum and minimum concentration of the chemical in the domain, i.e. $\max(c(t)) - \min(c(t))$, seen in Fig. \ref{fig:MaxMinChem}. This can give us insight into the density of the aggregations and the development of the chemical throughout the system. At short times (inset of Fig. \ref{fig:MaxMinChem}), this function grows linearly for all three deposition functions. After the nonlinear saturation takes place, the functions grow at different rates. We see that $f_1$ continues to grow and outpaces both $f_2$ and $f_3$, which implies that the aggregations for $f_1$ continue to get more densely populated over time. Both $f_2$ and $f_3$ grow at a much slower rate than $f_1$ and reach a quasi-equilibrium state for this max-min function.

To quantitatively show differences between the structure of these deposition functions, we utilize the radial distribution function, $g(r)$. The radial distribution function (RDF) measures the probability of finding a pair of plankton a distance of $r$ apart relative to the same probability for a completely uniform distribution of plankton. The RDF $g(r)$ is defined as (\citealt{Haile1993, Wang2019})
\begin{equation}\label{eq:RDF}
    g(r) = \frac{\ell^2}{N_p (N_p - 1)} \sum\limits_{i = 1}^{N_p}\sum\limits_{\substack{j=1 \\ j \neq i}}^{N_p} \delta\left(r - r_{ij} \right)
\end{equation}
Here, $N_p$ represents the total number of plankton in the system, $\ell$ represents the size of the domain, $\delta(x)$ represents the delta function, and $r_{ij} := \| \vec{x_i} - \vec{x_j} \|$ represents the Euclidean distance between plankton $i$ and $j$. If $g(r) = 1$, the probability is equivalent to a uniform distribution. If $g(r) > 1$, there is a greater probability of finding plankton a distance $r$ apart relative to a uniform distribution and thus, we expect an accumulation of plankton at a distance $r$ apart. Similarly, $g(r) < 1$ corresponds to a dearth of plankton a distance $r$ away. To calculate $g(r)$ numerically at a specific time $t$, we calculate all distances $r_{ij}$, bin these distances into a histogram, and normalize with respect to the total number of distances calculated. 

The RDFs are calculated for the various deposition functions over time and are shown in Fig. \ref{fig:RDFs}. At $t = 15$, we see that plankton become slightly more likely to find other organisms close to them $(r < 1)$ for all deposition function. However, the plot of $g(r)$ changes to reflect the dense aggregations formed by $f_1$ as seen in plots $t = 225$ and $t = 300$. We are much more likely to find plankton close together $(r < 3)$ for simulations with $f_1$ than for $f_2$ and $f_3$. We see that the graphs of $f_2$ and $f_3$ stagnate at around $t = 150$ whereas $f_1$ continues to evolve until the final timesteps. One can also estimate the average size of the aggregations in the system by looking at the first instance where $g(r) = 1$. We can then estimate at $t = 300$ that the aggregation diameter for $f_1$ is roughly $3.5$ units, whereas $f_2$ and $f_3$ are both around $2.7$ units, even though the aggregation for $f_2$ is more densely packed since $g(r)$ is greater. 
\begin{figure}[t!]
    \centering
    \subfloat[Total chemical, $C(t)$]{
        \includegraphics[width=.48\linewidth]{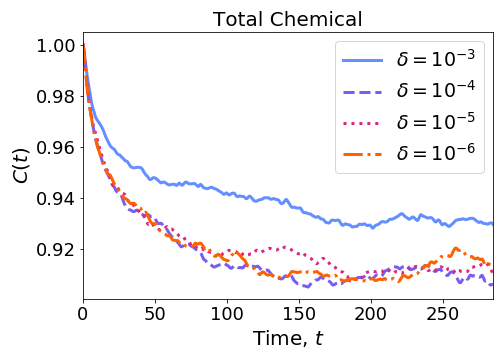}\label{fig:ChemicalTotDelt}}
        \hspace{.01\linewidth}
    \subfloat[Mean value of $|\nabla c|$ at a given $\rho$, $t = 300$]{
        \includegraphics[width=.48\linewidth]{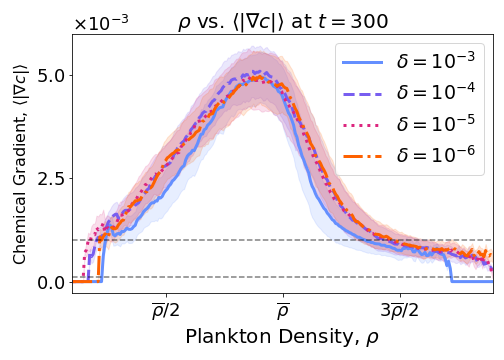}\label{fig:DeltGradc}}
    \caption{Plots \textbf{a} and \textbf{b} show characteristics of varied $\delta$. Plot \textbf{a} shows the the total chemical, $C(t)$, which is calculated in \eqref{eq:TotChemeq}. Plot \textbf{b} is the mean magnitude of the chemical gradient at a certain plankton density with error bars representing one standard deviation. The two dashed lines represent $10^{-4}$ and $10^{-3}$. Parameters for these simulations can be found in Appendix \ref{sec:2DParams}.}
    \label{fig:Deltas2D}
\end{figure}

Finally, we see that $\delta$ has a profound impact on the development of the chemical and plankton density. In Fig. \ref{fig:Deltas2D}, we use varied $\delta$ along with the switch deposition function $f_2$. In Fig. \ref{fig:ChemicalTotDelt}, we see that as $\delta \to 0$, the total chemical in the system decreases and for $\delta < 10^{-4}$, we see no significant difference in the function $C(t)$. Next, we recognize that tumbling is only relevant where plankton are present. We observe that larger aggregations tend to have near constant density interiors surrounded by transition areas leading to voids. Thus, we expect $\delta$ to play a role in these high density, non-constant interiors. To this end, Fig. \ref{fig:DeltGradc} shows the mean value of the magnitude of the chemical gradient, $\langle | \nabla c | \rangle$, at a given plankton density $\rho$ at $t = 300$. The gradient in the interiors are on the order of $10^{-3}$, and when $\delta \approx \langle | \nabla c | \rangle$, we expect the nonlinear dynamics to change as is shown in the plot. We infer that $\delta$ only has an impact inside the aggregation, where the gradient is low and the density is high. Therefore, the chemotactic sensitivity and numerical-regularization parameter $\delta$ alters both the linear stability of the steady state (see Section \ref{sec:Analytic2D}) and the nonlinear, saturated configuration of the plankton. 

\section{Discussion}\label{sec:Discussion}

We now discuss our results from the previous sections by connecting them towards the experimental observations in Section \ref{sec:Experiments}.

Our stability analysis in Figs. \ref{fig:RealRoot} and Figs. \ref{fig:RealRoot2D} gives us insight of how a constant solution $(\overline{c}, \overline{\rho})$ can become unstable and initiate plankton aggregations. If we decrease $d_1$ or increase $d_2$, we are able to make the system more unstable, i.e., as the chemical diffuses more slowly or decays more rapidly, the plankton are more likely to aggregate. Also, the importance of $\delta$ is evident, implying that the ability of the model organism to respond to weak gradients could be one of the keys to understanding the persistent patterns observed in nature, as shown by \cite{Sengupta2017}.

Turning towards the long term behavior seen in Sections \ref{sec:NumRes1D} and \ref{sec:NumRes2D}, by implementing only run-and-tumble and autochemotaxis, we gather rich behavior in one- and two-dimensions that is very similar to the the pattern formation we see in our model organism (see Fig. \ref{fig:Progression}). From simulations seen in Figs. \ref{fig:SSNOnConst} and \ref{fig:2DEV}, we observe several aggregations initially form and merge, but as aggregations coarsened and expanded, merging events occurred on much larger timescales. Specifically in the two-dimensional case, simulations such as Fig. \ref{fig:PlankDens}, which show the rapid organization of aggregations and subsequent long-time saturation towards larger, non-circular aggregations, reflect the structure and evolution patterns we see in experiments highlighted in Section \ref{sec:Experiments}. As expected, we are able to predict that the chemical concentration is highest near the peaks of the aggregations. However, the shape and evolution of these plankton aggregations and chemical concentrations are dependent on the properties of the deposition function and $\delta$.

We studied three different deposition functions and saw very different results for all three. By comparing deposition functions that admit similar most unstable wave numbers, we were able to focus on the non-linear development after the initial perturbation becomes apparent. In both the one- and two-dimensional cases, the constant function, $f_1$, is the most comparable to the evolution dynamics seen in our simulations when light saturated the environment. The switch function, $f_2$, forces the plankton to aggregate more slowly and merge on a different timescale than the other deposition functions. The linear switch function, $f_3$, produced  Even though $f_2$ and $f_3$ produced more unique aggregations, this may accurately reflect the speed of merging when there is no light as in Fig. \ref{fig:DarkExp}.

Lastly, our long time simulations varying $\delta$ in both the one- and two-dimensional cases stress the significance of a plankton's chemotactic sensitivity. As $\delta$ increases, the aggregations are much more spread out and less defined. As well, as $\delta$ decreases, plankton are more sensitive to the environment and thus can sustain aggregations mild gradients. Compare these results to experimental long term behavior as in Fig. \ref{fig:T60}. The aggregations after 60 minutes are relatively varied, as some are densely populated whereas others are spread out much more. This may imply that the plankton's chemotactic sensitivity is influenced by biological mechanics related to a plankton's response towards light. 

\section{Conclusions and Future Work}\label{sec:Conclusions}

We conducted several experiments to understand aggregation patterns of phytoplankton using our model organism \textit{Heterosigma akashiwo}. From these experiments, it is evident that chemotaxis plays an important role in plankton aggregations, as these aggregations form in absence of external fluid flow and light. To describe this system, we have constructed a one- and two-dimensional partial differential equation model of plankton with run-and-tumble motion and autochemotaxis. We explored the complex dynamics of the models by analytically describing the linear stability of constant steady states and numerically analyzing the nonlinear saturation and long time behavior. By varying parameters, we investigated the diverse dynamics of this model from stability analysis to structural differences at steady state. Since the deposition functions for plankton are unknown, we chose three simple deposition functions to simulate. Through simulations, these functions revealed striking differences in chemical and plankton aggregation patterns, some of which are similar to experimental observations. 

There is still much to be explored in this line of research. There are several significant biological and physical mechanisms, such as plankton photosynthesis, vertical migration towards and away from light, and a non-uniform turning kernel $T(\theta^{\prime}, \theta)$ which is biologically relevant, that need to be incorporated to understand plankton aggregation formation (\citealt{Dodson1997,Chen2018}). Extending this model to three dimensions may allow us to understand the extreme change in pattern formation seen in Fig. \ref{fig:DeepPlate}. We also will conduct more experimentation with collaborators to find parameters and deposition function which describes the behavior of \textit{Heterosigma akaswhio}.

\section*{Acknowledgments}

We would like to acknowledge Prof. Cathy Coyne and Monica Acerbi for providing us with {\em Heterosigma akashiwo} and food and advice for keeping them alive. We also acknowledge the Center for Computational Biology, directed by Dr. Michael Shelley, for their expertise in biological mechanisms underlying the plankton. As well, thank you to undergraduates Claire Lubash and Diana Li for their help in laboratory experimentation and numerical simulations.  

\section*{Declarations}
\subsection*{Conflicts of interest}
The authors declare that they have no conflict of interest.

\subsection*{Code Availability}
For the code which produced all figures and videos in this paper, see our \href{https://github.com/ricknussell31/AutochemotaxisFiles}{GitHub repository}.
For the videos and Online Resources, please go to this \href{https://drive.google.com/drive/folders/1AXWG54FCohTtUPOxIOTmcmddl9CVx9-F?usp=sharing}{folder}.

\bibliographystyle{spbasic} 
\bibliography{BibPlank}

\appendix 

\section{Laboratory experiments}
\label{sec:experiments}
 The \textit{Heterosigma akashiwo} were collected from the Delaware Bay in Lewes, Delaware and kept in beakers filled with seawater. The water was infused with nutrients with a combination of 1 mL of \ce{NaNO3}, 1 mL of \ce{NaH2PO4*H2O}, 1 mL of \ce{Na2SiO3*9H2O}, and .5 mL of f/2 trace metal and vitamin solution per 100 mL of sea water. All of the laboratory experiments were conducted in the MECLab in Ewing Hall at the University of Delaware. We used two different types of dishes for the experiments: a shallow dish and a deep dish. The dimensions of the shallow dish are $9.5 \times 9.5 \times 1$ centimeters and for the deep dish, the dimensions are $7 \times 7 \times 11$ centimeters. For the light source, we used an 80 watt incandescent bulb that attains approximately 1000 lumens, and we changed the angle and distance of the light source for the various experiments. The light source was 7-13 cm away from the dish and was at a 45$^{\circ}$-70$^{\circ}$ angle of incidence. To capture the movement of the plankton, we utilized an Allied Vision Mako G-30 camera equipped with Vimba Viewer software to take photos every 5 seconds over the duration of the experiments. Windows Media Player was used to stitch the photos together and construct the movies. For videos, visit the YouTube page for the \href{https://www.youtube.com/channel/UCheboMplIjxYTkjJuSBsBbA}{University of Delaware Math Plankton Team}.
 
 \section{Proofs of Propositions}\label{sec:Proofs}
\subsection{Proposition \ref{prop1}} \label{sec:Proof1}
\begin{proof}
To analyze the roots of $\eqref{eq:charact}$ as $k \to \infty$, we divide $\eqref{eq:charact}$ by $k^4$ and define $\alpha = k^{-1}$ to obtain
\begin{equation}\label{eq:charatdiv}
\begin{split}
        \alpha^4 \lambda^3 + \left[ d_1 \alpha^2 + (1 - d_3) \alpha^4 \right] \lambda^2 & + \left[ (d_1 + 1) \alpha^2 - d_3 \alpha^4 \right]\\
        & + \left( d_1 - \left( d_3 + \frac{\overline{c} d_2}{\delta} \right) \alpha^2 \right) = 0.
    \end{split}
\end{equation}
As $\alpha \to 0$ in $\eqref{eq:charatdiv}$, a singular perturbation occurs. Performing dominant balance on the first two terms, defining a variable $y = \alpha^{2} \lambda$, and substituting this $y$ representation into \eqref{eq:charatdiv}, we derive the cubic equation
\begin{equation}\label{eq:finalpert}
\begin{split}
         y^3 + \left[ d_1 + (1 - d_3) \alpha^2 \right] y^2 & + \left[ (d_1 + 1) \alpha^2 - d_3 \alpha^4 \right]\\
        & + \left( d_1 \alpha^2 - \left( d_3 + \frac{\overline{c} d_2}{\delta} \right) \alpha^4 \right) = 0.
        \end{split}
\end{equation}
We assume that $y$ is a solution of the form $y = \sum\limits_{n = 0}^{\infty} b_n \alpha^n$, where $b_n \in \mathbb{C}$. Substituting into $\eqref{eq:finalpert}$ and setting the terms corresponding to similar powers of $\alpha$ equal to 0, we obtain
\begin{align*}
    \alpha^0:&  \hspace{2mm} b_0^3 + d_1b_0^2 = 0 \Rightarrow b_0 = 0, 0, -d_1 \\
    \alpha^1:& \hspace{2mm} b_1(3b_0^2 + 2d_1 b_0) = 0 \\
    \alpha^2:& \hspace{2mm} b_0^2 \left( 3 b_2 - d_3 - 1\right) + b_0 \left( 1 + 3b_1 + d_1 + 2 d_1 b_2 \right) + d_1\left( 1 + b_1^2 \right) = 0 \\
    \alpha^3:& \hspace{2mm} b_1^3 + 3b_0^2 b_3 + 2 d_1 b_1b_2 + (d_1 + 1)b_1 + \left(6b_1b_2 + 2 d_1 b_3 - 2 (1 + d_3) b_1 \right) b_0 = 0.
\end{align*}
Taking $b_0 = 0$, the $\alpha^1$ equation gives no information about $b_1$. We then use the $\alpha^2$ equation to arrive at
\begin{equation}\label{eq:b1}
    d_1(1 + b_1^2) = 0 \Longrightarrow b_1 = -i, i.
\end{equation}
Taking $b_1 = - i$ and using the $\alpha^3$ equation, we get
\begin{equation}\label{eq:b2negi}
    i - 2id_1 b_2 - 2i d_1 - i = 0 \Longrightarrow b_2 = -\frac{1}{2}.
\end{equation}
Now taking $b_1 = i$ and using the $\alpha^3$ equation, we get
\begin{equation}
    -i + 2i d_1 b_2 + 2_i d_1 + i = 0 \Longrightarrow b_2 = -\frac{1}{2}.
\end{equation}
Therefore, the three roots are
\begin{align*}
       & y_1 = - d_1 + o(\alpha) && \Rightarrow  && \lambda_1 = -d_1 k^2 + \mathcal{O}(k) \\
       & y_2 = -i \alpha + -\frac{1}{2}\alpha^2 + o(\alpha^3) && \Rightarrow && \lambda_2 = - i k - \frac{1}{2} + \mathcal{O}\left(k^{-1}\right)\\
        & y_3 = i \alpha + -\frac{1}{2}\alpha^2 + o(\alpha^3) && \Rightarrow && \lambda_3 =  i k - \frac{1}{2} + \mathcal{O}\left(k^{-1}\right).
\end{align*}
Since $\text{Re}(\lambda_1) \to - \infty$, $\text{Re}(\lambda_2) \to -\frac{1}{2}$, and $\text{Re}(\lambda_3) \to -\frac{1}{2}$ as $k \to \infty$, then $k_u \to -\frac{1}{2}$. \qed
\end{proof}

\subsection{Proposition \ref{prop2}} \label{sec:Proof2}
\begin{proof}
\indent Recall that $R(k) = \max \{\text{Re}\left(\lambda_i\right)\}_{i =1}^3$. Setting $k = 0$ in \eqref{eq:charact}, we have 
\begin{equation}\label{eq:k0root}
    \lambda^3 + (1-d_3) \lambda^2 - d_3 \lambda = 0 \Longrightarrow \lambda = 0, -1, d_3
\end{equation}
Therefore, $R(0) = 0$ or $R(0) = d_3$, depending on the sign of $d_3$. In order to find the maximum of $R(k)$, we find the $k$ such that $\frac{\partial \lambda}{\partial k} = 0$. Using implicit differentiation on $\eqref{eq:charact}$, we obtain
\begin{dmath}\label{eq:firstderiv}
    \frac{\partial \lambda}{\partial k} \left( 3 \lambda^2 + 2 \lambda(d_1 k^2 - d_3 + 1) + (d_1 + 1)k^2 - d_3 \right) + \left(2 d_1 k \lambda^2 + 2k(d_1 + 1)\lambda + 4d_1k^3 - 2k\left(d_3 + \frac{\overline{c} d_2 }{\delta}\right) \right) = 0.
\end{dmath}
When $\frac{\partial \lambda}{\partial k} = 0$, we get 
\begin{equation}\label{eq:deriv0}
    2k\left( d_1 \lambda^2 + (d_1 + 1)\lambda + 2d_1k^2 - \left(d_3 + \frac{\overline{c} d_2}{\delta}\right)\right) =0.
\end{equation}
Since $k \geq 0$, \eqref{eq:deriv0} implies that $k = 0$ is a critical point and therefore, $R'(0) = 0$. To classify this critical point, we seek to find $R''(0)$. Solving for the second derivative at $k = 0$ using implicit differentiation of \eqref{eq:firstderiv}, we obtain 
\begin{equation}\label{eq:secondderiv0}
\frac{\partial^2 \lambda}{\partial k^2}\Big\vert_{k = 0} = -\frac{2 \left(d_3 + \frac{\overline{c} d_2}{\delta} - d_1 \lambda^2\right)}{3 \lambda^2 + 2 \lambda (1-d_3) - d_3}.
\end{equation}
If $R''(0) > 0$, then by Proposition \ref{prop1} and the fact that $R(k)$ is a continuous function, there must be a $k_m$ such that $R'(k_m) = 0$ and $R(k_m) > 0$. This will give us the bounds in \eqref{eq:ConditionProp2}.

To this end, we first assume that $d_3 > 0$. Using $\lambda = d_3$, \eqref{eq:secondderiv0} becomes
\begin{equation}\label{eq:secondderiv02}
\frac{\partial^2 \lambda}{\partial k^2}\Big\vert_{k = 0,\lambda = d_3} =  R''(0) = \frac{2 \left(d_3 + \frac{\overline{c} d_2}{\delta} - d_1 d_3^2\right)}{d_3^2 + d_3}.
\end{equation}
To have $R''(0) > 0$, the parameters must satisfy
\begin{equation}\label{eq:DDd3}
\left[ d_3 - \frac{1}{2 d_1} \left( 1 -  \sqrt{1+\frac{4 d_1 d_2 \overline{c}}{\delta}}\right)\right] \left[ d_3 - \frac{1}{2 d_1} \left( 1 + \sqrt{1+\frac{4 d_1 d_2 \overline{c}}{\delta}}\right)\right] < 0.
\end{equation}
Since $d_3 > 0$, \eqref{eq:DDd3} implies $d_3 < \frac{1}{2d_1} \left[1 + \sqrt{1 + \frac{4 \overline{c} d_1 d_2}{\delta}} \right]$, which gives the upper bound of \eqref{eq:ConditionProp2}. \\
\indent Now we assume that $d_3 <0$ and thus $R(0) = 0$. Substituting $\lambda = 0$ into \eqref{eq:secondderiv02}, we obtain
\begin{equation}\label{eq:secondderiv03}
\frac{\partial^2 \lambda}{\partial k^2}\Big\vert_{k = 0,\lambda = 0} =  R''(0) = \frac{2 \left(d_3 + \frac{\overline{c} d_2}{\delta} \right)}{-d_3}.
\end{equation}
This implies that $d_3 > -\frac{\overline{c} d_2}{\delta}$ will make $R''(0) > 0$, which gives the lower bound in \eqref{eq:ConditionProp2}. \qed 
\end{proof}

\subsection{Proposition \ref{prop:Prop3}} \label{sec:Proof3}
\begin{proof}
Consider the PDE for the chemical concentration, $c$, defined in \eqref{eq:2nond1} where $\rho(\vec{x}, t) = \sum\limits_{i = 1}^{N_p} \delta\left(\vec{x} - \vec{x}_i(t) \right)$. We solve this system on a periodic domain $D = [0, \ell] \times [0,\ell]$ and define $C(t)$ as the total chemical over time (see \eqref{eq:TotChemeq}). Integrating \eqref{eq:2nond1} over $D$, noting that $f(c) = f_1(c) = \gamma$, and using the steady state $d_2 \overline{c}  = f(\overline{c}) \overline{\rho}$, we obtain 
\begin{dgroup}
\begin{dmath}
\iint\limits_{D} c_t \, d \vec{x} = \iint\limits_{D} \left[ d_1 \Delta c - d_2 c + f(c) \rho \right] \, d \vec{x}
\end{dmath}
\begin{dmath}
C' =  d_2 C + \iint\limits_{D} f(c) \rho \, d \vec{x}
\end{dmath}
\begin{dmath}\label{eq:ODEPlankD}
C' -  d_2 C = \frac{\overline{\rho} \ell^2}{N_p}\sum\limits_{i = 1}^{N_p} f\left(c \left( \vec{x}_i(t)\right) \right) 
\end{dmath}
\begin{dmath}
C' -  d_2 C = \frac{\overline{\rho} \ell^2}{N_p} \left( \gamma N_p \right)
\end{dmath}
\begin{dmath}
C' -  d_2 C = \overline{\rho} \gamma \ell^2
\end{dmath}
\begin{dmath}
C(t) = \frac{\overline{\rho} \gamma \ell^2}{d_2} + A_1 e^{-d_2 t},
\end{dmath}
\begin{dmath}
C(t) = \frac{\overline{c} \gamma \ell^2}{f(\overline{c})} + A_1 e^{-d_2 t}
\end{dmath}
\begin{dmath}\label{eq:InitODE}
C(t) = \overline{c} \ell^2 + A_1 e^{-d_2 t}
\end{dmath}
\end{dgroup}
where $A_1$ is a constant and $\frac{\overline{\rho} \ell^2}{N_p}$ in \eqref{eq:ODEPlankD} is the plankton density coefficient described in \eqref{eq: PlankDensCoef}. Solving for $A_1$ using the initial condition of $C(0) = C_0$, the solution to the ODE becomes
\begin{dmath}\label{eq:FinalODE}
C(t) = \overline{c} \ell^2 + (C_0 - \overline{c} \ell^2)e^{-d_2 t}.
\end{dmath}
Therefore, as $t \to \infty$, $C(t) \to \overline{c} \ell^2$ and if $C_0 = \overline{c} \ell^2$, $C(t) = \overline{c} \ell^2$ for all $t$. \qed
\end{proof}
 \section{One Dimensional Numerical Simulations}
 \label{sec:1DNum}
 
 \subsection{Numerical Methods}\label{sec:1DMethods}
 
In this section, we discuss the details of the numerical scheme utilized in solving \eqref{eq:nd11}-\eqref{eq:nd12} on the interval $[0,\ell]$ with periodic boundary conditions

\begin{equation}\label{eq:BCs1D}
    \rho(0,t) = \rho(\ell, t), \quad \rho_x(0,t) = \rho_x(\ell,t), \quad c(0,t) = c(\ell,t), \quad c_x(0,t) = c_x(\ell,t)
\end{equation}
We define $\vec{D}$ as the $C_m \times C_m$ spectral Chebyshev differentiation matrix on the interval $[0,\ell]$ with periodic boundary conditions, where $C_m \in \mathbb{Z}^+$. We also denote $\boldsymbol{\rho}^n$ and $\vec{c}^n$ as the $n$th time step of the plankton density and chemical concentration, $\Delta t$ as the temporal mesh size, and $\Delta x = \frac{C_m}{\ell}$ as the spatial mesh size.
 
For \eqref{eq:nd11}, we utilize a pseudospectral method: a Crank-Nicolson method for time-stepping and a Chebyshev spectral method for the spatial mesh. To deal with the non-linearity of autochemotactic, we use a forward Euler method solely for that term. The scheme can be written as follows:
\begin{dmath*}
\frac{\vec{c}^{n+1} - \vec{c}^n}{\Delta t} = \frac{d_1}{2} \vec{D}^2\left[\vec{c}^n + \vec{c}^{n +1}\right] - \frac{d_2}{2} \left[\vec{c}^n + \vec{c}^{n + 1}\right] + f(\vec{c}^n) \boldsymbol{\rho}^n 
\end{dmath*}
\begin{dmath}\label{eq:1DCN}
\vec{c}^{n + 1} = \left[ \vec{I}  - (\Delta t/2) \left( d_1 \vec{D}^2 - d_2 \vec{I} \right) \right]^{-1} \left[ \left( \vec{I}  + (\Delta t/2) \left( d_1 \vec{D}^2 - d_2 \vec{I} \right) \right) \vec{c}^n + \Delta t f(\vec{c}^n) \boldsymbol{\rho}^n \right],
\end{dmath}
where $\vec{I}$ is the identity matrix of appropriate size. For \eqref{eq:nd12}, we utilize a pseudospectral method: a leap-frog method for the temporal mesh and a Chebyshev spectral method for the spatial mesh. The scheme can be written as follows:
\begin{dmath*}
\frac{\boldsymbol{\rho}^{n+1} - 2\boldsymbol{\rho}^n + \boldsymbol{\rho}^{n-1}}{(\Delta t)^2} + \frac{\boldsymbol{\rho}^{n+1} - \boldsymbol{\rho}^{n-1}}{2\Delta t} = \vec{D}^2 \boldsymbol{\rho}^{n} - \vec{D} \left( \dfrac{\vec{D} \vec{c}^n}{\sqrt{(\vec{D}\vec{c}^n)^2 + \delta^2}} \boldsymbol{\rho}^n \right)
\end{dmath*}
\begin{dmath}
\boldsymbol{\rho}^{n+1} = \frac{1}{1 + \Delta t/2} \left\lbrace 2 \boldsymbol{\rho}^n + \boldsymbol{\rho}^{n - 1}\left( \Delta t/2  - 1 \right) + (\Delta t)^2 \left[ \vec{D}^2 \boldsymbol{\rho}^n - \vec{D} \left( \dfrac{\vec{D} \vec{c}^n}{\sqrt{(\vec{D}\vec{c}^n)^2 + \delta^2}} \boldsymbol{\rho}^n \right) \right] \right\rbrace
\end{dmath}
To initialize the system, we define values for $d_1, d_2, \delta, \ell,\gamma, c_0, T_0, C_m, \Delta t$, and $f(c)$. We initialize $c$ and $\rho$ as constant solutions which satisfy the constant steady state solution of \eqref{eq:nd11}-\eqref{eq:nd12}, which is $d_2 \overline{c} = f(\overline{c}) \overline{\rho}$. For simulations in Section \ref{sec:NumRes1D}, we utilize $\ell = 5$, $\Delta t = 0.004$, and $C_m = 301$, which is spatially and temporally stable.  

\subsection{Parameters}\label{sec:1DParams}
Table \ref{tab:1DParams} shows all parameters used in the one dimensional simulations in Section \ref{sec:Analysis1D}. Details about a specific variable or notation can be found in Tables \ref{tab:ParametersD} and \ref{tab:ParametersND}.
\begin{table}[H]
\centering
\renewcommand\arraystretch{1.3}
\begin{tabular}{|c||c|c|c|c|c|c|c|c|c|}
\hline
Figure & $d_1$ & $d_2$ & $\ell$ & $\delta$ & $f(c)$ & $\gamma$ & $c_0$ & $T_0$ & $\overline{c}$\\
    \hline
    \hline
Fig. \ref{fig:MaxReald1} & Var. & 1 & -- & 0.01 & $f_2$ & 0.01 & 0.05 & 0.03 & 0.12\\
    \hline
Fig. \ref{fig:MaxReald2} & 1 & Var. & -- & 0.01 & $f_2$ & 0.01 & 0.05 & 0.03 & 0.12 \\
    \hline
Fig. \ref{fig:1DStabled1} & Var. & Var. & 6 & 0.01 & $f_2$ & 0.01 & 0.05 & 0.03 & 0.12\\
    \hline
Fig. \ref{fig:1DStabled2} & Var. & Var. & 6 & 0.012 & $f_2$ & 0.01 & 0.05 & 0.03 & 0.12 \\
    \hline
Fig. \ref{fig:SSNOnConst} & 0.2 & 4 & 5 & 0.015 & $f_1$ & 0.01 & -- & -- & 0.2\\
    \hline
Fig. \ref{fig:ConstantDep} & 0.2 & 4 & 5 & 0.015 & $f_1$ & 0.01 & -- & -- & 0.2 \\
    \hline
Fig. \ref{fig:SwitchDep} & 0.2 & 4 & 5 & 0.015 & $f_2$ & 0.01 & 0.25 & 0.04 & 0.2 \\
    \hline
Fig. \ref{fig:LinSwitchDep} & 0.2 & 4 & 5 & 0.015 & $f_3$ & 0.01 & 0.25 & 0.04 & 0.2 \\
    \hline
Fig. \ref{fig:DeltaPlank} & 1 & 2 & 5 & Var. & $f_1$ & 0.01 & -- & -- & 0.12 \\
    \hline
Fig. \ref{fig:DeltaChem} & 1 & 2 & 5 & Var. & $f_1$ & 0.01 & -- & -- & 0.12 \\
    \hline
\end{tabular}
\caption{Parameters used in the 1D analysis and simulations discussed in Section \ref{sec:Analysis1D}. ``Var." means the variable was varied and ``--" denotes unused variables.} 
\label{tab:1DParams}
\end{table}

 \section{Two Dimensional Numerical Simulations}
 \label{sec:2DNum}
 \subsection{Methods}\label{sec:2DMethods}
In this section, we describe the numerical scheme used to solve the two-dimensional model from the initial discussion in Section \ref{sec:NumRes2D}. Unlike the 1D scheme, we discretize the density field as a linear combination of moving particles representing an ensemble of plankton. Therefore, we introduce the plankton density coefficient, $\overline{\rho}_d$. Recall that $\overline{\rho}$ is the constant steady state solution to the 2D system, and define the particle density as $\frac{N_p}{\ell^2}$. In order to scale the autochemotactic term appropriately, each of the $N_p$ particles must represent a small aggregation of plankton. To represent the specified density on the computational domain using $N_p$ particles, we require that
\begin{equation}\label{eq: PlankDensCoef}
    \frac{\overline{\rho}}{\overline{\rho}_d} = \frac{N_p}{\ell^2} \Longrightarrow \overline{\rho}_d = \frac{\overline{\rho} \ell^2}{N_p}.
\end{equation}
Thus, for the autochemotactic term, each particle deposits chemical with strength $\overline{\rho}_d f(c)$. 

The chemical field is approximated on a mesh. W  we use a periodic $[0,\ell] \times [0, \ell]$ domain, and we define $\widetilde{\mathbf{D}}$ as the $2C_m \times 2C_m$ two-dimensional spectral Chebyshev differentiation matrix on the domain $[0,\ell]\times [0,\ell]$ with periodic boundary conditions. We denote $\boldsymbol{\rho}^n$ and $\vec{c}^n$ as the $n$th time step of the plankton density and chemical concentration, where $\boldsymbol{\rho}^n$ is a summation of point sources
\begin{equation}\label{eq:PlankDens}
    \boldsymbol{\rho}^n = \sum\limits_{i = 1}^{N_p} \delta(\vec{x}_i^n),
\end{equation}
and $\vec{x}_i^n = \langle x_i^n, y_i^n \rangle$ is the position of the $i$th plankton at the $n$th time step. We define $\Delta t$ as the temporal mesh size, $\Delta x = \Delta y = \frac{C_m}{\ell}$ as the spatial mesh size, and $\vec{x}_{\ell}$ and $\vec{y}_{\ell}$ as the spatial mesh for the chemical concentration, $c$, in the $x$ and $y$ directions, respectively.
 
To initialize the system, recall that the constant steady state solution to $\eqref{eq:2nond1}$-$\eqref{eq:2nond3}$ is $d_2 \overline{c} = f(\overline{c}) \overline{\rho}$. We select a $\overline{c}$, $d_2$, and $f(c)$ and then solve for $\overline{\rho} = \frac{d_2 \overline{c}}{f(\overline{c})}$. All $N$ organisms are randomly distributed throughout the domain and given a random orientation $\theta_i$.

For the $(n+1)$th iteration, we complete the following steps:
\begin{enumerate}
    \item Decide if plankton will tumble to a new direction using the probability described in $\eqref{eq:ProbND}$. If a plankton does tumble, their new direction is selected randomly from a uniform distribution $\theta_i^{n + 1} = [0,2\pi)$.
    \item Move all $N_p$ plankton to their new position using \eqref{eq:UpdatePos}, which updates $\rho^n$ to $\rho^{n + 1}$.
    \item Compute $S_i^{n+1} := f\left(c^{n}\left(\vec{x}_i^{n + 1}\right)\right)$, where $\vec{x}_i^{n + 1}$ is the $i$th plankton's position calculated from step 2. We compute $c^{n}\left(\vec{x}_i^{n + 1}\right)$, the chemical concentration at a particle position, using a bivariate spline interpolation.
    \item Deposit the chemical. In order to do deposit the chemical, we construct the matrix
    \begin{equation}\label{eq:AutoChemoTerm}
        \vec{S}^{n+1} =  \sum\limits_{i =1}^{N_p} \dfrac{S_i^{n + 1}}{4 \pi \sigma^2} \exp\left[ \frac{(x_i^{n +1} - \vec{x}_{\ell})^2 + (y_i^{n + 1} - \vec{y}_{\ell})^2}{4 \sigma^2}\right],
    \end{equation}
    where the variance is $\sigma^2 = K \Delta t d_1$ and $K$ is an suitable constant. In the simulations, we have chosen $K = 3$.
    \item Evolve the chemical concentration using a Crank-Nicholson method, similar to the one-dimensional case in \eqref{eq:1DCN}:
    \begin{equation}\label{eq:2DCN}
    \vec{c}^{n + 1} = \left[ \vec{I}  - \frac{\Delta t}{2} \left( d_1 \widetilde{\mathbf{D}}^2 - d_2 \vec{I} \right) \right]^{-1} \left[ \left( \vec{I}  + \frac{\Delta t}{2} \left( d_1 \widetilde{\mathbf{D}}^2 - d_2 \vec{I} \right) \right) \vec{c}^n + \Delta t \overline{\rho}_d \vec{S}^{n+1} \right],
    \end{equation}
    where $\vec{I}$ is the identity matrix of appropriate size.
\end{enumerate}

For all simulations in Section \ref{sec:Analysis2D}, we utilize $\Delta t = 0.2$, $\Delta x = \Delta y = 0.025$, $C_m = 200$, and $\ell = 10$. Using these parameters along with all others used in simulations described in Section \ref{sec:NumRes2D}, the largest wave number we are able to resolve has a magnitude of $|\vec{k}| \approx 41$. This mode, along with others with larger magnitude, are stable in all of the simulations, and thus all important wave numbers are resolved.

\subsection{Parameters}\label{sec:2DParams}
Table \ref{tab:2DParams} shows all parameters used in the two dimensional simulations in Section \ref{sec:Analysis2D}. Details about a specific variable or notation can be found in Tables \ref{tab:ParametersD} and \ref{tab:ParametersND}. For all simulations, the number of plankton used was $N_p = 1.6 \times 10^5$.
\begin{table}[H]
\centering
\renewcommand\arraystretch{1.3}
\begin{tabular}{|c||c|c|c|c|c|c|c|c|c|}
\hline
Figure & $d_1$ & $d_2$ & $\ell$ & $\delta$ & $f(c)$ & $\gamma$ & $c_0$ & $T_0$ & $\overline{c}$\\
    \hline
    \hline
Fig. \ref{fig:FirstUnstable2d} & Var. & 1 & -- & 0.001 & $f_2$ & 0.01 & 0.05 & 0.03 & 0.012 \\
    \hline
Fig. \ref{fig:SecondUnstable2d} & 1 & Var. & -- & 0.001 & $f_2$ & 0.01 & 0.05 & 0.03 & 0.012  \\
    \hline
Fig. \ref{fig:Ncheck} & 1 & 1 &  -- & 0.001 & $f_2$ & 0.01 & 0.05 & 0.03 & 0.012 \\
    \hline
Fig. \ref{fig:c0move} & 1 & 1 & -- & 0.001 & $f_1$/$f_2$ & 0.01 & Var. & 0.03 & 0.012  \\
    \hline
 Fig. \ref{fig:2DStability1delta} & Var. & Var.& 10 & 0.001 & $f_2$ & 0.01 & 0.05 & 0.03 & 0.012 \\
    \hline
Fig. \ref{fig:2DStability2delta} & Var. & Var. & 10 & 0.004 & $f_2$ & 0.01 & 0.05 & 0.03 & 0.012 \\
    \hline
Fig. \ref{fig:2DEV} & 0.1 & 4 & 10 & 0.001 & $f_1$ & 0.01 & -- & -- & 0.01\\
\hline
Figs. \ref{fig:PlankDens}-\ref{fig:MetricsDeps} & 0.1 & 4 & 10 & 0.001 & $f_1$ & 0.01 & -- & -- & 0.01 \\
\hline
Figs. \ref{fig:PlankDens}-\ref{fig:MetricsDeps} & 0.1 & 4 & 10 & 0.001 & $f_2$ & 0.01 & 0.012 & 0.0007 & 0.01\\
\hline 
Figs. \ref{fig:PlankDens}-\ref{fig:MetricsDeps} & 0.1 & 4 & 10 & 0.001 & $f_3$ & 0.01 & 0.012 & 0.008 & 0.01\\
\hline 
Fig. \ref{fig:Deltas2D} & 0.1 & 4 & 10 & Var. & $f_2$ & 0.01 & 0.012 & 0.0007 & 0.01\\
\hline 
\end{tabular}
\caption{Parameters used in the 2D analysis and simulations discussed in Section \ref{sec:Analysis2D}. ``Var." means the variable was varied and ``--" denotes unused variables.} 
\label{tab:2DParams}
\end{table}

\end{document}


\begin{flushleft}
\textbf{\href{https://drive.google.com/file/d/1Z1c6ByRpX89Ey64_CbMkfnkweeG_CiZ8/view?usp=sharing}{Online Resource 1:}} Simulation of the one-dimensional system (7a)-(7b) from a constant steady state using the constant deposition function $f(c) = f_1(c)$. The red line denotes $\rho$, the plankton density, and the blue line denotes $c$, the chemical concentration. To see parameters used for this simulation, see Appendix C.2 and Table 3. \\
\vspace{5mm}

\textbf{\href{https://drive.google.com/file/d/1ElmZNh7AUS3BuLLMaiPPbSz9t0v2ui5O/view?usp=sharing}{Online Resource 2:}} Simulation of the one-dimensional system (7a)-(7b) from a constant steady state using the switch deposition function $f(c) = f_2(c)$. The red line denotes $\rho$, the plankton density, and the blue line denotes $c$, the chemical concentration. To see parameters used for this simulation, see Appendix C.2 and Table 3. \\
\vspace{5mm}

\textbf{\href{https://drive.google.com/file/d/109TwWza5H7JFiJv9fzEFMLF2CzXMupA-/view?usp=sharing}{Online Resource 3:}} Simulation of the one-dimensional system (7a)-(7b) from a constant steady state using the constant deposition function $f(c) = f_1(c)$. The red line denotes $\rho$, the plankton density, and the blue line denotes $c$, the chemical concentration. To see parameters used for this simulation, see Appendix C.2 and Table 3. \\
\vspace{5mm}

\textbf{\href{https://drive.google.com/file/d/1BdfVCqi99g94p3rJiiW8RuO9KAeLgWqK/view?usp=sharing}{Online Resource 4:}} Plankton density of the two-dimensional system described in Section 5.2 from a constant steady state using the constant deposition function $f(c) = f_1(c)$. To see parameters used for this simulation along with the numerical techniques used, see Appendix D and Table 4.\\
\vspace{5mm} 

\textbf{\href{https://drive.google.com/file/d/1q79_w1pHBOEle5FMKjc0saa1d74uZsSt/view?usp=sharing}{Online Resource 5:}} Plankton density of the two-dimensional system described in Section 5.2 from a constant steady state using the switch deposition function $f(c) = f_2(c)$. To see parameters used for this simulation along with the numerical techniques used, see Appendix D and Table 4. \\
\vspace{5mm}

\textbf{\href{https://drive.google.com/file/d/1yBs5yypQQokwiTQFtFuv-l8mAD9D-X2S/view?usp=sharing}{Online Resource 6:}} Plankton density of the two-dimensional system described in Section 5.2 from a constant steady state using the linear switch deposition function $f(c) = f_3(c)$. To see parameters used for this simulation along with the numerical techniques used, see Appendix D and Table 4. \vspace{5mm}

\textbf{Contact Information:}
\begin{itemize}
    \item Article: Analysis and Simulation of a Novel Run-and-Tumble Model with Autochemotaxis
    \item Journal: Journal of Mathematical Biology
    \item Author Names: Nicholas J. Russell and Louis F. Rossi
    \item Corresponding Author: Nicholas J. Russell, Department of Mathematical Sciences, University of Delaware. Email: nrussell@udel.edu
    \end{itemize}
\end{flushleft}